\documentclass{article}
\usepackage[utf8]{inputenc}
\usepackage[spanish,english]{babel}
\usepackage{amsmath}
\usepackage{amsthm}
\usepackage{amssymb}
\usepackage{mathtools}
\usepackage{mathrsfs}
\usepackage[inline]{enumitem}
\usepackage{csquotes}
\usepackage{tikz-cd}
\usepackage{placeins}
\usepackage{xkeyval}
\usepackage{multicol}
\usepackage{calc}
\usepackage{subcaption}
\usepackage{algorithm}
\usepackage{algpseudocode}
\usepackage{array}
\usepackage[hyphens]{url}
\usepackage[style=apa]{biblatex}

\newcommand*\Let[2]{\State #1 $\gets$ #2}
\newcommand*\Forr[2]{\For{$#1 \textbf{ to } #2$}}
\algrenewcommand\alglinenumber[1]{{\sffamily\footnotesize#1:}}
\algrenewcommand\algorithmicrequire{\textbf{Precondition:}}
\algrenewcommand\algorithmicensure{\textbf{Postcondition:}}
\algrenewcommand\algorithmiccomment[1]{%
  \quad{\sffamily\footnotesize \# #1}}
\setlist[enumerate,1]{%
  label={\normalfont(\alph*)},
  ref={\normalfont(\alph*)}
}

\DeclareRobustCommand
  \vdotss{\vbox{\baselineskip3pt\lineskiplimit0pt
    \kern6pt\hbox{.}\hbox{.}\hbox{.}}}

\newcommand\restr[2]{{
  \left.\kern-\nulldelimiterspace 
  #1 
  \vphantom{\big|} 
  \right|_{#2} 
  }}

\DeclareMathOperator{\Ima}{Im}
\DeclareMathOperator{\rank}{rank}
\DeclareMathOperator{\Codim}{codim}
\DeclareMathOperator{\Hom}{Hom}

\newcommand\dime[1][V]{\underline{\dim}\,#1}
\DeclareMathOperator{\End}{End}

\newcommand\carc[1]{\mathcal{#1}}

\newcommand\recmatrix[6][]{%
\begin{tikzpicture}[baseline=(current bounding box.center),remember picture]
\matrix[mimat,#1] (#6)
{
#2 \& #3 \\
#4 \& #5 \\
};
\draw (#6-1-1.north west) -- (#6-2-1.south west);
\draw (#6-1-1.north east) -- (#6-2-1.south east);
\draw (#6-1-2.north east) -- (#6-2-2.south east);
\draw (#6-1-1.north west) -- (#6-1-2.north east);
\draw (#6-1-1.south west) -- (#6-1-2.south east);
\draw (#6-2-1.south west) -- (#6-2-2.south east);
\end{tikzpicture}%
}

\newcommand\knob[1]{$k$\nobreakdash-\hspace{0pt}#1}

\newcommand\geomatrix[4]{$%
  \setlength\arraycolsep{4pt}
  \begin{array}{
  |>{\centering\arraybackslash}p{8pt}
  |>{\centering\arraybackslash}p{8pt}|
  }
  \hline
  #1 & #2
  \\\hline
  #3 & #4
  \\\hline
  \end{array}$}
  
\newcommand\tttmatrix[9][]{%
\begin{tikzpicture}[#1]
\matrix[mimat,inner sep=0pt,text width=3em] (mim)
{
#2 \& #3 \\
#4 \& #5 \\
};
\draw (mim-1-1.north west) -- (mim-2-1.south west);
\draw (mim-1-1.north east) -- (mim-2-1.south east);
\draw (mim-1-2.north east) -- (mim-2-2.south east);
\draw (mim-1-1.north west) -- (mim-1-2.north east);
\draw (mim-1-1.south west) -- (mim-1-2.south east);
\draw (mim-2-1.south west) -- (mim-2-2.south east);

\node[anchor=north,font=\footnotesize] 
  at ([yshift=-0.5ex]mim.south) 
  {$n\geq #6$};
\node[anchor=south,font=\footnotesize] 
  at ([yshift=0.5ex]mim.north) (end) 
  {$\mathcal{E}\simeq #9$};
\node[anchor=south,font=\footnotesize] 
  at ([yshift=1ex]end) (dim) 
  {$q(d)=#7$};
\node[anchor=south,font=\footnotesize] 
  at ([yshift=1ex]dim) 
  {$d=#8$};
\end{tikzpicture}%
}
\newcommand\tttmatrixwide[9][]{%
\begin{tikzpicture}[#1]
\matrix[mimat,inner sep=0pt,text width=4.5em] (mim)
{
#2 \& #3 \\
#4 \& #5 \\
};
\draw (mim-1-1.north west) -- (mim-2-1.south west);
\draw (mim-1-1.north east) -- (mim-2-1.south east);
\draw (mim-1-2.north east) -- (mim-2-2.south east);
\draw (mim-1-1.north west) -- (mim-1-2.north east);
\draw (mim-1-1.south west) -- (mim-1-2.south east);
\draw (mim-2-1.south west) -- (mim-2-2.south east);

\node[anchor=north,font=\footnotesize] 
  at ([yshift=-0.5ex]mim.south) 
  {$n\geq #6$};
\node[anchor=south,font=\footnotesize] 
  at ([yshift=0.5ex]mim.north) (end) 
  {$\mathcal{E}\simeq #9$};
\node[anchor=south,font=\footnotesize] 
  at ([yshift=1ex]end) (dim) 
  {$q(d)=#7$};
\node[anchor=south,font=\footnotesize] 
  at ([yshift=1ex]dim) 
  {$d=#8$};
\end{tikzpicture}%
}

\newcommand\sqmatrix[4]{%
\begin{tikzpicture}[baseline=(current bounding box.center)]
\matrix[
    matrix of math nodes,
    ampersand replacement=\&,
    text width=1.5em,
    text height=1em,
    text depth=0.3em,
    align=center
]
(mim)
{
#1 \& #2 \\
#3 \& #4 \\
};
\draw (mim-1-1.north west) -- (mim-2-1.south west);
\draw (mim-1-1.north east) -- (mim-2-1.south east);
\draw (mim-1-2.north east) -- (mim-2-2.south east);
\draw (mim-1-1.north west) -- (mim-1-2.north east);
\draw (mim-1-1.south west) -- (mim-1-2.south east);
\draw (mim-2-1.south west) -- (mim-2-2.south east);
\end{tikzpicture}%
}

\newcommand\hormatrix[3][]{%
\begin{tikzpicture}[#1]
\matrix[mimat,inner sep=0pt] (mimat)
{
#2 \& #3 \\
};
\draw (mimat-1-1.north west) -- (mimat-1-1.south west);
\draw (mimat-1-1.north east) -- (mimat-1-1.south east);
\draw (mimat-1-2.north east) -- (mimat-1-2.south east);
\draw (mimat-1-1.north west) -- (mimat-1-2.north east);
\draw (mimat-1-1.south west) -- (mimat-1-2.south east);

\end{tikzpicture}%
}

\newtheorem{theo}{Theorem}[section]
\newtheorem{prop}[theo]{Proposition}
\newtheorem{coro}[theo]{Corollary}

\theoremstyle{definition}

\theoremstyle{remark}
\newtheorem{rema}{Remark}

\usetikzlibrary{
  calc,
  matrix,
  positioning,
  decorations.pathreplacing,
  decorations.markings,
  arrows,
  arrows.meta,
  shapes,
  backgrounds,
  babel,
  quotes
}

\tikzset{
  copo/.style={draw,circle,inner sep=1.5pt},
  texto/.style={draw=none,inner sep=2pt},
  mimatriz/.style={
    anchor=base,
    align=center,
    text depth=0.8ex,
    text height=2.2ex,
    text width=2em,
  matrix of math nodes,
  nodes in empty cells
  },
  mimat/.style={
    matrix of math nodes,
    nodes in empty cells,
    ampersand replacement=\&,
    text width=3.4em,
    text height=1.3em,
    text depth=0.6em,
    align=center
  },
  corto/.style={shorten >= 3pt,shorten <= 3pt},
  myarrow/.style={
    postaction={%
    decorate,
    decoration={markings,
    mark=at position #1 with {\arrow{>}}}}
  }
}

\numberwithin{equation}{section}

\newcommand\ile[1][n]{I^{\leftarrow}_{#1}}
\newcommand\iri[1][n]{I^{\rightarrow}_{#1}}
\newcommand\ido[1][n]{I^{\downarrow}_{#1}}
\newcommand\iup[1][n]{I^{\uparrow}_{#1}}

\DeclareMathOperator{\rep}{rep}
\newcommand\repSk{\rep(\carc{S},k)}
\newcommand\repS{\rep\carc{S}}

\newcommand\inte[1]{\{1,\dots,#1\}}

\title{On a classification problem for a quiver of type $\widetilde{A}_{3}$}
\author{Ivon Dorado\thanks{Universidad Nacional de Colombia, sede Bogotá, \texttt{iadoradoc@unal.edu.co}} \and Gonzalo Medina\thanks{Universidad Nacional de Colombia, sede Manizales, \texttt{gmedinaar@unal.edu.co}}}

\begin{document}

\maketitle
\begin{abstract}
We present a new solution to the classification problem for the category of representations of a quiver of type $\widetilde{A}_{3}$. Our approach uses linear algebra techniques which lead us to a reduction that allows to use induction. As an application, the solution to the classical Kronecker problem and its contragredient version are obtained in an elementary way. We also describe the endomorphism rings for the indecomposable representations and an algorithm that shows how to reconstruct their matrix form from some graphic invariants.\\[1ex]
\textbf{Keywords:} Indecomposable representation; Quiver of type $\widetilde{A}_{3}$; Matrix presentation; Endomorphism ring.
\end{abstract}

\section*{Introduction}
In this paper, we present a new solution to the classification problem of four linear operators defined between four finite-dimensional vector spaces over a field. This problem is equivalent to finding all the indecomposable representations, up to isomorphism, for the following quiver $\carc{S}$ of type $\widetilde{A}_{3}$ (recall that $\widetilde{A}_{n}$ denotes the affine or extended Dynkin diagram with $n+1$ vertices):
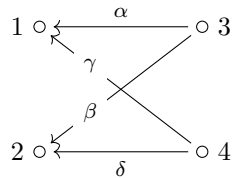
\[
\begin{tikzpicture}[node distance=1.5cm and 2cm,baseline=(current bounding box.center)]
\node[copo,label={left:$1$}] (1) {};
\node[copo,below=of 1,label={left:$2$}] (2) {};
\node[copo,right=of 1,label={right:$3$}] (3) {};
\node[copo,below=of 3,label={right:$4$}] (4) {};
\draw[->,shorten >= 3pt,shorten <= 3pt] 
  (3) -- 
  node[above,font=\footnotesize] {$\alpha$} 
  (1);
\draw[->,shorten >= 3pt,shorten <= 3pt] 
  (3) -- 
  node[fill=white,pos=0.7] {\phantom{$\beta$}}
  node[font=\footnotesize,pos=0.7] {$\beta$}
  (2);
\draw[->,shorten >= 3pt,shorten <= 3pt] 
  (4) -- 
  node[fill=white,pos=0.7] {\phantom{$\gamma$}}
  node[font=\footnotesize,pos=0.7] {$\gamma$}
  (1);
\draw[->,shorten >= 3pt,shorten <= 3pt] 
  (4) -- 
  node[below,font=\footnotesize] {$\delta$}
  (2);
\end{tikzpicture}
\]
The problem was first solved in~\textcite{Naz1967} using a matrix-based approach.

We use a matrix approach and a \textquote{reduction} mechanism that will allow us to use induction. As an application, we show how to easily obtain the solution to some classical classification problems. We also obtain the endomorphism rings for the indecomposable representations.

We will use some well-known facts from the theory of representation of quivers and from standard linear algebra, but to make the exposition as self-contained as possible, we include in the text all necessary definitions. 

The paper is organized as follows. In Section~\ref{sec:indecomposable} we state our main result about the indecomposable representations and we also give some preliminaries and introduce the notation that will be used throughout the paper. In Section~\ref{sec:solution} we present the solution to the classification problem for the indecomposable representations and introduce some graphical invariants for some of the indecomposable objects. In Section~\ref{sec:endomorphisms} we obtain our second main result: the classification of the endomorphism rings associated to the indecomposable representations. In Section~\ref{sec:subproblems}, we show how to use our solution to easily obtain solutions to the Kronecker problem and its contragredient versions. Finally, in the appendix, we present an algorithm that shows that the graphical invariants introduced in Section~\ref{sec:solution} are indeed enough to reconstruct the corresponding indecomposable representations. 

\section{Indecomposable representations for a quiver of type $\widetilde{A}_{3}$}
\label{sec:indecomposable}

For a given field $k$, \textit{representations of $\carc{S}$} are $8$\nobreakdash-tuples having the form
\[
V=(V_{1},V_{2},V_{3},V_{4},f_{\alpha},f_{\beta},f_{\gamma},f_{\delta}),
\] 
where $V_{1}$, $V_{2}$, $V_{3}$, $V_{4}$ are vector spaces over the field $k$, and $f_{\alpha}\colon V_{3}\to V_{1}$, $f_{\beta}\colon V_{3}\to V_{2}$, $f_{\gamma}\colon V_{4}\to V_{1}$, $f_{\delta}\colon V_{4}\to V_{2}$ are \knob{linear} maps. A representation is \textit{finite-dimensional} if all four $V_{1}$, $V_{2}$, $V_{3}$, and $V_{4}$ are finite dimensional vector spaces over $k$. For such a representation, its \textit{dimension vector $d$} is the element of $\mathbb{Z}^{4}$ given by
\[
d=\dime =(d_{1},d_{2},d_{3},d_{4}),
\]
where $d_{i}=\dim_{k}V_{i}$, for all $i\in\inte{4}$. The \textit{Tits quadratic form} of $\carc{S}$ is the quadratic form $q=q_{\carc{S}}\colon\mathbb{Z}^{4}\to \mathbb{Z}$ given by
\[
q(x_{1},\ldots,x_{4})=
\sum_{i=1}^{4}x_{i}^{2}
  -x_{3}x_{1}-x_{3}x_{2}-x_{4}x_{1}-x_{4}x_{2}.
\]
An element $(x_{1},\ldots,x_{4})\in\mathbb{Z}^{4}-\{0\}$ is a \textit{real root} of the form $q$ if $q(x_{1},\ldots,x_{4})=1$, and it is an \textit{imaginary root} of $q$ if $q(x_{1},\ldots,x_{4})=0$.

In the following, we will refer to arbitrary representations of $\carc{S}$ given by:
\begin{gather*}
V=(V_{1},V_{2},V_{3},V_{4},
  f_{\alpha},f_{\beta},f_{\gamma},f_{\delta}) \\
\shortintertext{and} W=(W_{1},W_{2},W_{3},W_{4},
  g_{\alpha},g_{\beta},g_{\gamma},g_{\delta})
\end{gather*}
A \textit{morphism} $l\colon V\to W$ is a collection $l=(l_{1},\ldots,l_{4})$  of four \knob{linear} maps $l_{i}\colon V_{i}\to W_{i}$, for $i\in\inte{4}$, such that the following diagram commutes:
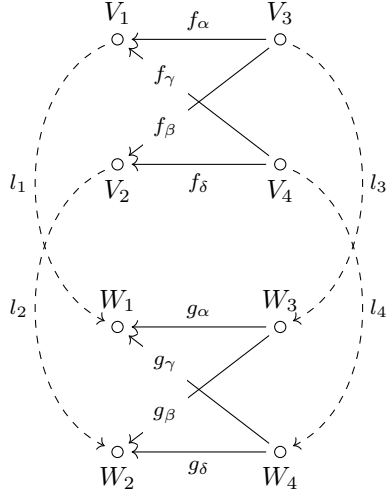
\[
\begin{tikzpicture}[node distance=1.5cm and 2cm]
\node[copo,label={above:$V_{1}$}] (v1) {};
\node[copo,below=of v1,label={below:$V_{2}$}] (v2) {};
\node[copo,right=of v1,label={above:$V_{3}$}] (v3) {};
\node[copo,below=of 3,label={below:$V_{4}$}] (v4) {};
\draw[->,shorten >= 3pt,shorten <= 3pt] 
  (v3) -- 
  node[above,font=\footnotesize] {$f_{\alpha}$} 
  (v1);
\draw[->,shorten >= 3pt,shorten <= 3pt] 
  (v3) -- 
  node[fill=white,pos=0.72] {\phantom{$f_{\beta}$}}
  node[font=\footnotesize,pos=0.72] {$f_{\beta}$}
  (v2);
\draw[->,shorten >= 3pt,shorten <= 3pt] 
  (v4) -- 
  node[fill=white,pos=0.72] {\phantom{$f_{\gamma}$}}
  node[font=\footnotesize,pos=0.72] {$f_{\gamma}$}
  (v1);
\draw[->,shorten >= 3pt,shorten <= 3pt] 
  (v4) -- 
  node[below,font=\footnotesize] {$f_{\delta}$}
  (v2);

\node[copo,below=2cm of v2,label={above:$W_{1}$}] (w1) {};
\node[copo,below=of w1,label={below:$W_{2}$}] (w2) {};
\node[copo,right=of w1,label={above:$W_{3}$}] (w3) {};
\node[copo,below=of w3,label={below:$W_{4}$}] (w4) {};
\draw[->,shorten >= 3pt,shorten <= 3pt] 
  (w3) -- 
  node[above,font=\footnotesize] {$g_{\alpha}$} 
  (w1);
\draw[->,shorten >= 3pt,shorten <= 3pt] 
  (w3) -- 
  node[fill=white,pos=0.72] {\phantom{$g_{\beta}$}}
  node[font=\footnotesize,pos=0.72] {$g_{\beta}$}
  (w2);
\draw[->,shorten >= 3pt,shorten <= 3pt] 
  (w4) -- 
  node[fill=white,pos=0.72] {\phantom{$g_{\gamma}$}}
  node[font=\footnotesize,pos=0.72] {$g_{\gamma}$}
  (w1);
\draw[->,shorten >= 3pt,shorten <= 3pt] 
  (w4) -- 
  node[below,font=\footnotesize] {$g_{\delta}$}
  (w2);

\foreach \index/\direct in {1,2}
{
\draw[->,shorten >= 3pt,shorten <= 3pt,dashed]
  (v\index) to[out=205,in=155]
  node[left,font=\footnotesize] {$l_{\index}$}
  (w\index);  
}
\foreach \index/\direct in {3,4}
{
\draw[->,shorten >= 3pt,shorten <= 3pt,dashed]
  (v\index) to[out=-25,in=25]
  node[right,font=\footnotesize] {$l_{\index}$}
  (w\index);  
}
\end{tikzpicture}
\]
The morphism $l=(l_{1},\ldots,l_{4})$ is an \textit{isomorphism} if all four operators $l_{i}$ are linear isomorphisms; if this is the case, we say that $V$ and $W$ are \textit{isomorphic} and write $V\simeq W$.

A direct sum can be defined naturally between representations in the following way:
\begin{align*}
V\oplus W
&=(V_{1}\oplus W_{1},V_{2}\oplus W_{2},V_{3}\oplus W_{3},V_{4}\oplus W_{4},\\
&\qquad f_{\alpha}\oplus g_{\alpha},f_{\beta}\oplus g_{\beta},f_{\gamma}\oplus g_{\gamma},f_{\delta}\oplus g_{\delta}).
\end{align*} 
A non-zero representation $V$ is \textit{indecomposable} if $V\simeq V'\oplus V''$ implies $V'=0$ or $V''=0$, where $0$ denotes the zero representation.
 
The representations of $\carc{S}$ and their morphisms form the category $\repSk$, or simply $\repS$, of representations of $\carc{S}$ over the field $k$. This category is abelian. For a proof, see~\textcite{AssSimSko2006-1}.
 
It will be useful for us to consider a natural duality for representations of $\carc{S}$. To each object $V$, we will associate its \textit{dual representation} 
\[
V^{\ast}=(V^{\ast}_{1},V^{\ast}_{2},V^{\ast}_{3},V^{\ast}_{4},f^{\ast}_{\alpha},f^{\ast}_{\beta},f^{\ast}_{\gamma},f^{\ast}_{\delta}),
\] 
where, for each $i\in\inte{4}$, the space $V^{\ast}_{i}$ is the algebraic dual of $V_{i}$, i.e., $V^{\ast}_{i}=\Hom_{k}(V_{i},k)$, and $f^{\ast}_{\alpha}$, $f^{\ast}_{\beta}$, $f^{\ast}_{\gamma}$, and $f^{\ast}_{\delta}$ are the corresponding dual operators. Notice that $V^{\ast}$ is a representation of the \textit{opposite quiver} $\carc{S}^{op}$:
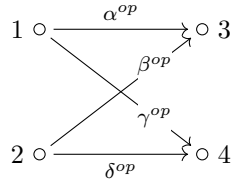
\[
\begin{tikzpicture}[node distance=1.5cm and 2cm,baseline=(current bounding box.center)]
\node[copo,label={left:$1$}] (1) {};
\node[copo,below=of 1,label={left:$2$}] (2) {};
\node[copo,right=of 1,label={right:$3$}] (3) {};
\node[copo,below=of 3,label={right:$4$}] (4) {};
\draw[<-,shorten >= 3pt,shorten <= 3pt] 
  (3) -- 
  node[above,font=\footnotesize] {$\alpha^{op}$} 
  (1);
\draw[<-,shorten >= 3pt,shorten <= 3pt] 
  (3) -- 
  node[fill=white,pos=0.28] {\phantom{$\beta$}}
  node[font=\footnotesize,pos=0.28] {$\beta^{op}$}
  (2);
\draw[<-,shorten >= 3pt,shorten <= 3pt] 
  (4) -- 
  node[fill=white,pos=0.28] {\phantom{$\gamma$}}
  node[font=\footnotesize,pos=0.28] {$\gamma^{op}$}
  (1);
\draw[<-,shorten >= 3pt,shorten <= 3pt] 
  (4) -- 
  node[below,font=\footnotesize] {$\delta^{op}$}
  (2);
\end{tikzpicture}
\]
obtained by keeping the same vertices as in $\carc{S}$ but reversing the direction of the arrows of $\carc{S}$.
 
We will present finite dimensional representations in matrix form: given a representation $V$, we first choose ordered bases for the spaces $V_{1}$, $V_{2}$, $V_{3}$, and $V_{4}$ and construct a block matrix of the form
\[
M=M_{V}=
\begin{tikzpicture}[baseline=(current bounding box.center)]
\matrix[matrix of math nodes,
    ampersand replacement=\&,
    text width=2.5em,
    text height=1.5em,
    text depth=0.8em,
    align=center
]
(mim)
{
M_{11} \& M_{12} \\
M_{21} \& M_{22} \\
};
\draw (mim-1-1.north west) -- (mim-2-1.south west);
\draw (mim-1-1.north east) -- (mim-2-1.south east);
\draw (mim-1-2.north east) -- (mim-2-2.south east);
\draw (mim-1-1.north west) -- (mim-1-2.north east);
\draw (mim-1-1.south west) -- (mim-1-2.south east);
\draw (mim-2-1.south west) -- (mim-2-2.south east);
\draw[decoration={brace,raise=5pt},decorate]
  (mim-1-2.north east) --
    node[xshift=7pt,anchor=west] {\footnotesize$\dim_{k}V_{1}$}
  (mim-1-2.south east);
\draw[decoration={brace,raise=5pt},decorate]
  (mim-2-2.north east) --
    node[xshift=7pt,anchor=west] {\footnotesize$\dim_{k}V_{2}$}
  (mim-2-2.south east);
\draw[decoration={brace,mirror,raise=5pt},decorate]
  (mim-2-1.south west) --
    node[yshift=-7pt,anchor=north] {\footnotesize$\dim_{k}V_{3}$}
  (mim-2-1.south east);
\draw[decoration={brace,mirror,raise=5pt},decorate]
  (mim-2-2.south west) --
    node[yshift=-7pt,anchor=north] {\footnotesize$\dim_{k}V_{4}$}
  (mim-2-2.south east);
\end{tikzpicture},
\] 
where the blocks $M_{11}$, $M_{21}$, $M_{12}$, and $M_{22}$ correspond to the matrices of the operators $f_{\alpha}$, $f_{\beta}$, $f_{\gamma}$, and $f_{\delta}$, respectively, with respect to the chosen bases. A matrix like $M_{V}$ before will be called a \textit{matrix presentation} for the representation $V$. Changing the chosen bases will now be modeled, in matrix form, by the following \textit{admissible transformations}:
\begin{enumerate}[noitemsep,label={\textsc{at}\arabic*.},ref={\textsc{at}\arabic*}]
\item\label{at1} Elementary row operations within each of the horizontal stripes of $M$.
\item\label{at2} Elementary column operations within each of the vertical stripes of $M$.
\end{enumerate}
Two representations $V$ and $W$ are isomorphic if any two of their matrix presentations $M_{V}$ and $M_{W}$ are \textit{equivalent}, denoted by $M_{V}\simeq M_{W}$, in the sense that one of them can be transformed into the other by applying a finite number of admissible transformations.

Direct sums of representations can be expressed in matrix form. Given two representations $V$ and $W$, with matrix presentations
\[
M_{V}=
\begin{tikzpicture}[baseline=(current bounding box.center)]
\matrix[matrix of math nodes,
    ampersand replacement=\&,
    text width=2em,
    text height=1.3em,
    text depth=0.6em,
    align=center
]
(mim)
{
M_{11} \& M_{12} \\
M_{21} \& M_{22} \\
};
\draw (mim-1-1.north west) -- (mim-2-1.south west);
\draw (mim-1-1.north east) -- (mim-2-1.south east);
\draw (mim-1-2.north east) -- (mim-2-2.south east);
\draw (mim-1-1.north west) -- (mim-1-2.north east);
\draw (mim-1-1.south west) -- (mim-1-2.south east);
\draw (mim-2-1.south west) -- (mim-2-2.south east);
\end{tikzpicture}
\qquad\text{and}\qquad
N_{W}=
\begin{tikzpicture}[baseline=(current bounding box.center)]
\matrix[matrix of math nodes,
    ampersand replacement=\&,
    text width=2em,
    text height=1.3em,
    text depth=0.6em,
    align=center
]
(mim)
{
N_{11} \& N_{12} \\
N_{21} \& N_{22} \\
};
\draw (mim-1-1.north west) -- (mim-2-1.south west);
\draw (mim-1-1.north east) -- (mim-2-1.south east);
\draw (mim-1-2.north east) -- (mim-2-2.south east);
\draw (mim-1-1.north west) -- (mim-1-2.north east);
\draw (mim-1-1.south west) -- (mim-1-2.south east);
\draw (mim-2-1.south west) -- (mim-2-2.south east);
\end{tikzpicture},
\]
their direct sum $V\oplus W$ has as one of its matrix presentations the following matrix, in which $0$ denotes zero blocks of the appropiate size:
\[
P_{V\oplus W}=M_{V} \oplus N_{W}=
\begin{tikzpicture}[baseline=(current bounding box.center)]
\matrix[matrix of math nodes,
    ampersand replacement=\&,
    text width=2em,
    text height=1.3em,
    text depth=0.6em,
    align=center
]
(mim)
{
M_{11} \& 0 \& M_{12} \& 0 \\
0 \& N_{11} \& 0 \& N_{12} \\
M_{21} \& 0 \& M_{22} \& 0 \\
0 \& N_{21} \& 0 \& N_{22} \\
};
\draw (mim-1-1.north west) -- (mim-4-1.south west);
\draw (mim-1-2.north east) -- (mim-4-2.south east);
\draw (mim-1-4.north east) -- (mim-4-4.south east);

\draw[dashed] (mim-1-1.north east) -- (mim-4-1.south east);
\draw[dashed] (mim-1-3.north east) -- (mim-4-3.south east);

\draw (mim-1-1.north west) -- (mim-1-4.north east);
\draw (mim-2-1.south west) -- (mim-2-4.south east);
\draw (mim-4-1.south west) -- (mim-4-4.south east);

\draw[dashed] (mim-1-1.south west) -- (mim-1-4.south east);
\draw[dashed] (mim-3-1.south west) -- (mim-3-4.south east);
\end{tikzpicture},
\]
We have everything necessary to state our main result:
\begin{theo}
\label{the:indecomposable2x2}
All indecomposable representations of $\carc{S}$ are exhausted, up to duality and up to permutations of horizontal and vertical stripes, by the matrix presentations listed in Figure~\ref{fig:indecomposableA3}.  
\end{theo}
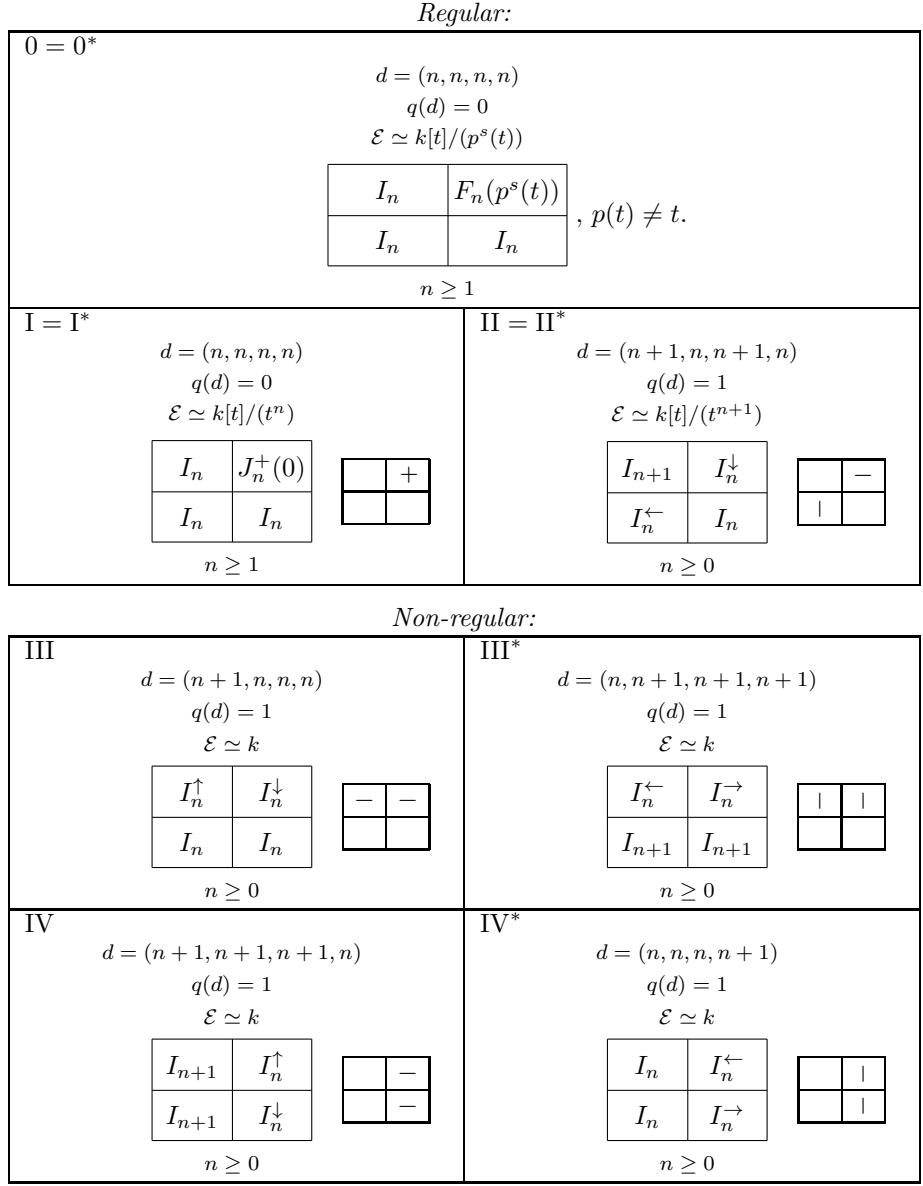
\begin{figure}
\begin{tabular}{%
  |>{\centering\arraybackslash}p{%
    \dimexpr0.5\textwidth-2\tabcolsep-1.5\arrayrulewidth\relax}|
  >{\centering\arraybackslash}p{%
    \dimexpr0.5\textwidth-2\tabcolsep-1.5\arrayrulewidth\relax}|
}
\multicolumn{2}{c}{\textit{Regular:}} \\
\hline
\multicolumn{2}{|>{\centering\arraybackslash}p{%
    \dimexpr\textwidth-4\tabcolsep-3\arrayrulewidth\relax}|}{%
    \rlap{$\mathrm{0}=\mathrm{0}^{\ast}$}\par
    \tttmatrixwide[remember picture,baseline=(current bounding box.center)]{I_{n}}{F_{n}(p^{s}(t))}{I_{n}}{I_{n}}{1}{0}{(n,n,n,n)}{k[t]/(p^{s}(t))}}%
\begin{tikzpicture}[remember picture,overlay]
\node[align=center,anchor=west,inner sep=0pt] at ([xshift=3pt]mim.east) 
{\text{, $p(t)\neq t$.}};
\end{tikzpicture}%
\\
\hline
\rlap{$\mathrm{I}=\mathrm{I}^{\ast}$}\par
\tttmatrix[baseline=(current bounding box.center),remember picture]{I_{n}}{J_{n}^{+}(0)}{I_{n}}{I_{n}}{1}{0}{(n,n,n,n)}{k[t]/(t^{n})}
\tikz[remember picture,overlay]{
\node[anchor=west] at ([xshift=8pt]mim.east) 
  {\geomatrix{}{\hspace{-3pt}{$-$}\hspace*{\widthof{$-$}*\real{-0.5}}\makebox[0pt][c]{\rotatebox[origin=c]{90}{$-$}}}{}{}};
}
&
\rlap{$\mathrm{II}=\mathrm{II}^{\ast}$}\par
\tttmatrix[baseline=(current bounding box.center),remember picture]{I_{n+1}}{\ido}{\ile}{I_{n}}{0}{1}{(n+1,n,n+1,n)}{k[t]/(t^{n+1})}
\tikz[remember picture,overlay]{
\node[anchor=west] at ([xshift=8pt]mim.east) 
  {\geomatrix{}{$-$}{\rotatebox[origin=c]{90}{$-$}}{}};
}
\\
\hline
\end{tabular}\par\medskip
\begin{tabular}{%
  |>{\centering\arraybackslash}p{%
    \dimexpr0.5\textwidth-2\tabcolsep-1.5\arrayrulewidth\relax}|
  >{\centering\arraybackslash}p{%
    \dimexpr0.5\textwidth-2\tabcolsep-1.5\arrayrulewidth\relax}|
}
\multicolumn{2}{c}{\textit{Non-regular:}} \\
\hline
\rlap{$\mathrm{III}$}\par
\tttmatrix[baseline=(current bounding box.center),remember picture]{\iup}{\ido}{I_{n}}{I_{n}}{0}{1}{(n+1,n,n,n)}{k}
\tikz[remember picture,overlay]{
\node[anchor=west] at ([xshift=8pt]mim.east) 
  {\geomatrix{$-$}{$-$}{}{}};
}
&
\rlap{$\mathrm{III}^{\ast}$}\par
\tttmatrix[baseline=(current bounding box.center),remember picture]{\ile}{\iri}{I_{n+1}}{I_{n+1}}{0}{1}{(n,n+1,n+1,n+1)}{k}
\tikz[remember picture,overlay]{
\node[anchor=west] at ([xshift=8pt]mim.east) 
  {\geomatrix{\rotatebox[origin=c]{90}{$-$}}{\rotatebox[origin=c]{90}{$-$}}{}{}};
}
\\
\hline
\rlap{$\mathrm{IV}$}\par
\tttmatrix[baseline=(current bounding box.center),remember picture]{I_{n+1}}{\iup}{I_{n+1}}{\ido}{0}{1}{(n+1,n+1,n+1,n)}{k}
\tikz[remember picture,overlay]{
\node[anchor=west] at ([xshift=8pt]mim.east) 
  {\geomatrix{}{$-$}{}{$-$}};
}
&
\rlap{$\mathrm{IV}^{\ast}$}\par
\tttmatrix[baseline=(current bounding box.center),remember picture]{I_{n}}{\ile}{I_{n}}{\iri}{0}{1}{(n,n,n,n+1)}{k}
\tikz[remember picture,overlay]{
\node[anchor=west] at ([xshift=8pt]mim.east) 
  {\geomatrix{}{\rotatebox[origin=c]{90}{$-$}}{}{\rotatebox[origin=c]{90}{$-$}}};
}
\\
\hline
\end{tabular}
\caption{Indecomposable representations of the quiver $\carc{S}$.}
\label{fig:indecomposableA3}. 
\end{figure}
In the next section, we will present the proof of Theorem~\ref{the:indecomposable2x2}. In the meantime, we will introduce some notations and definitions used in Figure~\ref{fig:indecomposableA3} and throughout this paper. 

Recall that using elementary row and column operations, any matrix $M$ of size $m\times n$ can be transformed into its standard form
\[
S(M,r)=
\begin{tikzpicture}[baseline=(current bounding box.center)]
\matrix[matrix of math nodes,
    ampersand replacement=\&,
    text width=1em,
    text height=1em,
    text depth=0.2em,
    align=center
]
(mim)
{
I_{r} \& 0 \\
0 \& 0 \\
};
\draw (mim-1-1.north west) -- (mim-2-1.south west);
\draw (mim-1-2.north east) -- (mim-2-2.south east);

\draw (mim-1-1.north west) -- (mim-1-2.north east);
\draw (mim-2-1.south west) -- (mim-2-2.south east);

\draw[dashed] (mim-1-1.north east) -- (mim-2-1.south east);
\draw[dashed] (mim-1-1.south west) -- (mim-1-2.south east);

\node[anchor=west,outer sep=0pt,inner sep=0pt] at (mim.east) {\small, $r\geq 1$,};
\end{tikzpicture}
\] 
where $r$ is the rank of $M$ (some blocks in the matrix above can be void).
 
For an integer $n\geq 1$, the matrices $\iup$, $\ido$ are the $(n+1)\times n$ matrices obtained by adjoining a row of zeros above, below, the identity matrix $I_{n}$; i.e.,
\[
\iup=
\begin{bmatrix}
0 & 0 & \dots & 0 \\
1 & 0 & \dots & 0 \\
0 & 1 & \dots & 0 \\
\vdots & \vdots & \ddots & \vdots \\
0 & 0 & \dots & 1
\end{bmatrix}_{(n+1)\times n}\qquad\text{and}\qquad
\ido=
\begin{bmatrix}
1 & 0 & \dots & 0 \\
0 & 1 & \dots & 0 \\
\vdots & \vdots & \ddots & \vdots \\
0 & 0 & \dots & 1 \\
0 & 0 & \dots & 0 \\
\end{bmatrix}_{(n+1)\times n}
\hspace{0pt-\widthof{$\scriptstyle(n+1)\times n$}}
\makebox[0pt][c]{.}
\hspace{\widthof{$\scriptstyle(n+1)\times n$}}
\]
Analogously, the matrices $\iri$, $\ile$ are the $n\times (n+1)$ matrices obtained by adjoining a column of zeros to the right, to the left, of the identity matrix $I_{n}$; i.e., 
\[
\iri=
\begin{bmatrix}
1 & 0 & \dots & 0 & 0 \\
0 & 1 & \dots & 0 & 0 \\
\vdots & \vdots & \ddots & \vdots & \vdots \\
0 & 0 & \dots & 1 & 0
\end{bmatrix}_{n\times (n+1)}\quad\text{and}\quad
\ile=
\begin{bmatrix}
0 & 1 & 0 & \dots & 0 \\
0 & 0 & 1 & \dots & 0 \\
\vdots & \vdots & \vdots & \ddots & \vdots \\
0 & 0 & 0 & \dots & 1
\end{bmatrix}_{n\times (n+1)}
\hspace{0pt-\widthof{$\scriptstyle n\times (n+1)$}}
\makebox[0pt][c]{.}
\hspace{\widthof{$\scriptstyle n\times (n+1)$}}
\]
For $n=0$, the matrices $\iup[0]$ and $\ido[0]$ are equal and they are \textquote{formal} matrices having one row and zero columns, and representing the linear operator $0\to k$. The matrices $\iri[0]$ and $\ile[0]$ are also equal and are \textquote{formal} matrices having zero rows and one column, and representing the linear operator $k\to 0$.

By $F_n(p^{s}(t))$ we will denote the Frobenius cell, also called rational canonical form cell, of order $n$ having the minimal polynomial $p^{s}(t)$, where $p(t)$ is monic and irreducible. In other words, $F_n(p^{s}(t))$ is the companion matrix of $p^{s}(t)$ (in particular, $n=s\cdot\deg p(t)$).

We denote by $J^{+}_{n}(0)$ or $J^{-}_{n}(0)$, respectively, the Jordan block of order $n$ with  eigenvalue $0$ and entries $1$ above or below the main diagonal, respectively.

In Figure~\ref{fig:indecomposableA3}, each type of indecomposable representation is accompanied by its corresponding dimension vector $d=\dime$ and by the value $q(d)$ of its Tits quadratic form. Each indecomposable object $V$ of one of the types
$\mathrm{II}$, $\mathrm{III}$, $\mathrm{IV}$, or of their dual types, is uniquely determined, up to isomorphism, by its dimension vector $d$, and each indecomposable of the types $\mathrm{0}$ or $\mathrm{I}$ is determined by the pair $(d,p(t))$ or $(d,t)$, respectively. Each indecomposable representation is also accompanied by its endomorphism ring $\mathcal{E}$ (see Section~\ref{sec:endomorphisms}).

In all matrices that have Jordan blocks, the block $J^{+}_{r}(0)$ can be replaced by $J^{-}_{r}(0)$ and the corresponding presentations are isomorphic.

With the exception of type $\mathrm{0}$, all other types of matrix presentations in~Figure~\ref{fig:indecomposableA3} include companion diagrams containing two symbols, either $-$ or $\rotatebox[origin=c]{90}{$-$}$. These graphical diagrams, along with the value $n$, entirely characterize the indecomposable objects in the sense that even without knowing the form of the matrix presentation, it can be constructed using the procedure presented in Algorithm~\ref{alg:graphicinvariant} in the Appendix. 

Notice that interchanging the symbols $-$ and $\rotatebox[origin=c]{90}{$-$}$ in a diagram produces the diagram for the corresponding dual presentation, up to permutations of vertical and horizontal stripes. 

\section{The solution to the classification problem for $\carc{S}$}
\label{sec:solution}
In this section we will prove Theorem~\ref{the:indecomposable2x2}, but first we need to establish an auxiliary result that will be essential in the proof.
\begin{prop}
\label{pro:reducedindecomposable}
Let $V = (V_{1},V_{2},V_{3},V_{4},f_{\alpha},f_{\beta},f_{\gamma},f_{\delta})$ be a representation of $\carc{S}$ such that $\dim_{k}(\Ima f_{\gamma})<\dim_{k}(V_{1})$, and let us define a new representation 
\[
V' =(V'_{1},V'_{2},V'_{3},V'_{4},
  f'_{\alpha},f'_{\beta},f'_{\gamma},f'_{\delta})
\]   
in the following way. Take $V'_{1}=\Ima f_{\gamma}$, and write $V_{1}=V'_{1}\oplus F$, for some complement $F$, with $\dim_{k}F=r\geq 1$. Let us take $B=\{u_{1},\ldots,u_{r}\}$ a set of $r$ linear independent vectors of $V_{3}$ such that $\langle f_{\alpha}(B)\rangle=\langle\lbrace f_{\alpha}(u_{1}),\ldots,f_{\alpha}(u_{r})\rbrace\rangle=F$, and extend this set to a basis $B\cup B'$ for $V_{3}$. Write $V_{3}=Y\oplus Z$, where $Y=\langle B\rangle$ and $Z=\langle B' \rangle$. Take $V'_{2}=V_{2}$, $V'_{3}=Z$, and $V'_{4}=V_{4}$. Linear operators $f'_{\alpha}$ and $f'_{\beta}$ are the restrictions given by $f'_{\alpha}=\restr{f_{\alpha}}{Z}$ and $f'_{\beta}=\restr{f_{\beta}}{Z}$. Furthermore, $f'_{\gamma}$ is $f_{\gamma}$ with its codomain changed to $V'_{1}$ (so $f'_{\gamma}=f_{\gamma}$) and $f'_{\delta}=f_{\delta}$. Then $V$ is indecomposable if and only if $V'$ is indecomposable.
\end{prop}
Notice that in our result we are considering a representation $V$ in which the spaces $V_{1}$ and $V_{3}$ are decomposed as in the diagram on the left. Then we state that $V$ is indecomposable if and only if the representation $V'$ on the right is.
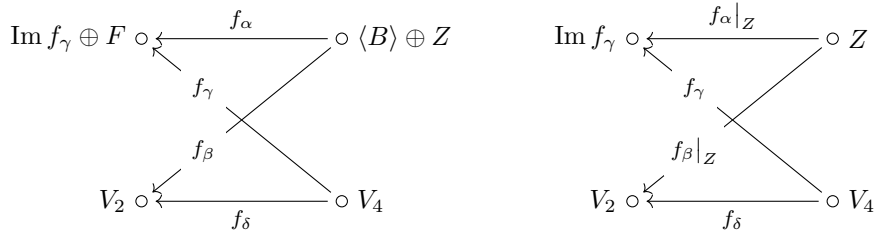
\[
\begin{tikzpicture}[node distance=2cm and 2.5cm,baseline=(current bounding box.center)]
\node[copo,label={left:$\Ima f_{\gamma}\oplus F$}] (1) {};
\node[copo,below=of 1,label={left:$V_{2}$}] (2) {};
\node[copo,right=of 1,label={right:$\langle B\rangle \oplus Z$}] (3) {};
\node[copo,below=of 3,label={right:$V_{4}$}] (4) {};
\draw[->,shorten >= 3pt,shorten <= 3pt] 
  (3) -- 
  node[above,font=\footnotesize] {$f_{\alpha}$} 
  (1);
\draw[->,shorten >= 3pt,shorten <= 3pt] 
  (3) -- 
  node[fill=white,pos=0.7] {\phantom{$f_{\beta}$}}
  node[font=\footnotesize,pos=0.7] {$f_{\beta}$}
  (2);
\draw[->,shorten >= 3pt,shorten <= 3pt] 
  (4) -- 
  node[fill=white,pos=0.7] {\phantom{$f_{\gamma}$}}
  node[font=\footnotesize,pos=0.7] {$f_{\gamma}$}
  (1);
\draw[->,shorten >= 3pt,shorten <= 3pt] 
  (4) -- 
  node[below,font=\footnotesize] {$f_{\delta}$}
  (2);

\node[copo,right=3.7cm of 3,label={left:$\Ima f_{\gamma}$}] (R1) {};
\node[copo,below=of R1,label={left:$V_{2}$}] (R2) {};
\node[copo,right=of R1,label={right:$Z$}] (R3) {};
\node[copo,below=of R3,label={right:$V_{4}$}] (R4) {};
\draw[->,shorten >= 3pt,shorten <= 3pt] 
  (R3) -- 
  node[above,font=\footnotesize] {$\restr{f_{\alpha}}{Z}$} 
  (R1);
\draw[->,shorten >= 3pt,shorten <= 3pt] 
  (R3) -- 
  node[fill=white,pos=0.7] {\phantom{$\restr{f}{Z}$}}
  node[font=\footnotesize,pos=0.7] {$\restr{f_{\beta}}{Z}$}
  (R2);
\draw[->,shorten >= 3pt,shorten <= 3pt] 
  (R4) -- 
  node[fill=white,pos=0.7] {\phantom{$f_{\gamma}$}}
  node[font=\footnotesize,pos=0.7] {$f_{\gamma}$}
  (R1);
\draw[->,shorten >= 3pt,shorten <= 3pt] 
  (R4) -- 
  node[below,font=\footnotesize] {$f_{\delta}$}
  (R2);
\end{tikzpicture}
\]
\begin{proof}
The verification that $V'$, as introduced above, is a well-defined representation of $\carc{S}$ is left to the reader.

$\Rightarrow$) If we assume that $V'$ is decomposable, say $V'\simeq M\oplus N$ with $M$ and $N$ non-zero representations
\begin{gather*}
M = 
(M_{1},M_{2},M_{3},M_{4},g_{\alpha},g_{\beta},g_{\gamma},g_{\delta}) \\
\shortintertext{and}
N = 
(N_{1},N_{2},N_{3},N_{4},h_{\alpha},h_{\beta},h_{\gamma},h_{\delta}),
\end{gather*}
then $V$ would also be since, by our construction, we have
\begin{align*}
V_{1} &= V'_{1}\oplus F = M_{1}\oplus N_{1}\oplus F, &
  V_{2} &= V'_{2} = M_{2}\oplus N_{2}, \\
V_{3} &= Y\oplus V'_{3} = Y\oplus M_{3}\oplus N_{3}, &
  V_{4} &= V'_{4}= M_{4}\oplus N_{4}.    
\end{align*}
As for the morphisms, we have the following:
\begin{align*}
f_{\alpha} &= 
  \bigl(\restr{f_{\alpha}}{Y}\bigr)\oplus \bigl(\restr{f_{\alpha}}{Z}\bigr) 
  = \bigl(\restr{f_{\alpha}}{Y}\bigr)\oplus f'_{\alpha} 
  = \bigl(\restr{f_{\alpha}}{Y}\bigr)\oplus g_{\alpha} \oplus h_{\alpha},\\
f_{\beta} &= 
  \bigl(\restr{f_{\beta}}{Y}\bigr)\oplus \bigl(\restr{f_{\beta}}{Z}\bigr) 
  = \bigl(\restr{f_{\beta}}{Y}\bigr)\oplus f'_{\beta} 
  = \bigl(\restr{f_{\beta}}{Y}\bigr)\oplus g_{\beta} \oplus h_{\beta},\\
f_{\gamma} &=
  f'_{\gamma} = g_{\gamma} \oplus h_{\gamma}, \\
f_{\delta} &=
  f'_{\delta} = g_{\delta} \oplus h_{\delta}.
\end{align*}
From here, we get
\[
V\simeq M\oplus 
  (N_{1}\oplus F,N_{2},N_{3}\oplus Y,N_{4},h_{\alpha}\oplus \restr{f_{\alpha}}{Y},h_{\beta}\oplus \restr{f_{\beta}}{Y},h_{\gamma},h_{\delta}),
\]
which shows that $V$ decomposes in a non-trivial way.

$\Leftarrow$) Conversely, if $V\simeq U\oplus W$, for non-zero representations $U$ and $W$ of $\carc{S}$, then we will show that $V'\simeq  U'\oplus W'$. Let us take
\begin{gather*}
U
=    
  (U_{1},U_{2},U_{3},U_{4},g_{\alpha},g_{\beta},g_{\gamma},g_{\delta})\\
\shortintertext{and}   
W
=    
  (W_{1},W_{2},W_{3},W_{4},h_{\alpha},h_{\beta},h_{\gamma},h_{\delta}).
\end{gather*}
In the following, we will use subscripts, such as in $F_{V}$, $F_{U}$, and $F_{W}$ to distinguish the different complements chosen in the construction of $V'$, $U'$, and $W'$, respectively.

Let us first prove that $V'_{1}=U'_{1}\oplus W'_{1}$. We have 
\[
V'_{1} 
  = \Ima f_{\gamma}
  = \Ima (g_{\gamma}\oplus h_{\gamma}) 
  = \Ima g_{\gamma} \oplus  \Ima h_{\gamma}
  = U'_{1}\oplus W'_{1}.  
\]
Moreover, 
\[
V'_{1}\oplus F_{V} 
= V_{1}
\simeq U_{1} \oplus W_{1}
\simeq (U'_{1}\oplus W'_{1})\oplus (F_{U}\oplus F_{W})
= V'_{1}\oplus (F_{U}\oplus F_{W}),
\]
so $F_{V}\simeq F_{U}\oplus F_{W}$. Let us now take linearly independent subsets $B\subseteq U_{3}$, and $C\subseteq W_{3}$ such that $\langle g_{\alpha}(B) \rangle=F_{U}$, and $\langle h_{\alpha}(C) \rangle=F_{W}$. Extend $B$ and $C$ to obtain $B\cup B'$ and $C\cup C'$, bases for $U_{3}$ and $W_{3}$, respectively. Take 
\[
Y_{U} = \langle B \rangle,\quad 
Z_{U} = \langle B' \rangle = U'_{3},\quad 
Y_{W} = \langle C \rangle,\quad\text{and}\quad 
Z_{W} = \langle C' \rangle = W'_{3}.
\] 
We have $\langle g_{\alpha}(B)\cup h_{\alpha}(C)\rangle=F_{U}\oplus F_{W}\simeq F_{V}$, so $V'_{3}=Y_{V}=Y_{U}\oplus Y_{W}=U'_{3}\oplus W'_{3}$.
Finally, since the second and fourth subspaces of the representations do not change under our construction, we immediately get $V'_{2} = V_{2} = U_{2} \oplus W_{2} = U'_{2} \oplus W'_{2}$  and $V'_{4}=V_{4}=U_{4}\oplus W_{4}=U'_{4}\oplus W'_{4}$. Let us now proceed to examine the morphisms.

Since $f_{\alpha}=g_{\alpha}\oplus h_{\alpha}$ and $V'_{1}=U'_{1}\oplus W'_{1}$, we obtain
\[
f'_{\alpha} = \restr{f_{\alpha}}{V'_{1}}
  =\restr{(g_{\alpha}\oplus h_{\alpha})}{U'_{1}\oplus W'_{1}}
  =\restr{g_{\alpha}}{U'_{1}}\oplus \restr{h_{\alpha}}{W'_{1}}
  =g'_{\alpha}\oplus h'_{\alpha}.
\]
The verification $f'_{\beta}=g'_{\beta}\oplus h'_{\beta}$ is completely analogous. To conclude our proof, note that according to our construction, $f_{\gamma}=f'_{\gamma}$, $g_{\gamma}=g'_{\gamma}$, and $h_{\gamma}=h'_{\gamma}$, so $f_{\gamma}=g_{\gamma}\oplus h_{\gamma}$ immediately yields $f'_{\gamma}=g'_{\gamma}\oplus h'_{\gamma}$ and a similar situation applies to the $\delta$ morphisms. We have thus shown that $V'=U'\oplus W'$ and this decomposition is not trivial since, otherwise, the summands $\begin{array}{|c|}\hline 1\\\hline 1\\\hline\end{array}$ would split off from $V$.
\end{proof}
Now we have all the necessary elements to prove our main result: 
\begin{proof}[\normalfont\bfseries Proof of Theorem~\ref{the:indecomposable2x2}]
Let $V=(V_{1},V_{2},V_{3},V_{4},f_{\alpha},f_{\beta},f_{\gamma},f_{\delta})$ be an indecomposable representation with matrix presentation 
\[
M_{V}=
\begin{tikzpicture}[baseline=(current bounding box.center)]
\matrix[matrix of math nodes,
    ampersand replacement=\&,
    text width=2em,
    text height=1.3em,
    text depth=0.6em,
    align=center
]
(mim)
{
M_{11} \& M_{12} \\
M_{21} \& M_{22} \\
};
\draw (mim-1-1.north west) -- (mim-2-1.south west);
\draw (mim-1-1.north east) -- (mim-2-1.south east);
\draw (mim-1-2.north east) -- (mim-2-2.south east);
\draw (mim-1-1.north west) -- (mim-1-2.north east);
\draw (mim-1-1.south west) -- (mim-1-2.south east);
\draw (mim-2-1.south west) -- (mim-2-2.south east);
\end{tikzpicture}.
\] 
Two cases are possible:
\begin{enumerate}[%
  label={\makebox[42pt][l]{Case~\arabic*:~}},
  ref={\arabic*},
  labelsep=0pt,
  leftmargin=*,
  itemindent=42pt
]
\item\label{ite:caseA} All four operators $f_{\alpha}$, $f_{\beta}$, $f_{\gamma}$, and $f_{\delta}$ are isomorphisms. In this case, the four blocks $M_{11}$, $M_{12}$, $M_{21}$, $M_{22}$ correspond to square non-singular matrices and, applying appropriate admissible transformations, three of those four blocks can be transformed into identity blocks. For definitiveness, let us say that $M_{11}$, $M_{21}$, and $M_{22}$ are the three said blocks. We have that $M_{V}$ has been reduced to a matrix of the form
\[
M_{V} \simeq
\begin{tikzpicture}[baseline=(current bounding box.center)]
\matrix[matrix of math nodes,
    ampersand replacement=\&,
    text width=1.5em,
    text height=1em,
    text depth=0.3em,
    align=center
]
(mim)
{
I_{n} \& S \\
I_{n} \& I_{n} \\
};
\draw (mim-1-1.north west) -- (mim-2-1.south west);
\draw (mim-1-1.north east) -- (mim-2-1.south east);
\draw (mim-1-2.north east) -- (mim-2-2.south east);
\draw (mim-1-1.north west) -- (mim-1-2.north east);
\draw (mim-1-1.south west) -- (mim-1-2.south east);
\draw (mim-2-1.south west) -- (mim-2-2.south east);
\end{tikzpicture}.
\] 
In the following diagram, we examine which admissible transformations can be applied to the block $S$, while preserving the identity blocks. Elementary row (or column) operations on a horizontal (or vertical) stripe correspond to premultiplication (or postmultiplication) by a non-singular matrix.
\[
\begin{tikzpicture}[inner sep=0pt]
\node (a) 
{\recmatrix{I_{n}}{S}{I_{n}}{I_{n}}{J1}};
\node[right=of a] (b) 
{\recmatrix{X}{XS}{I_{n}}{I_{n}}{J2}};
\node[right=of b] (c) 
{\recmatrix{I_{n}}{XS}{X^{-1}}{I_{n}}{J3}};
\node[below=of c] (d) 
{\recmatrix{I_{n}}{XS}{I_{n}}{X}{J4}};
\node[left=of d] (e) 
{\recmatrix[text width=3.5em]{I_{n}}{XSX^{-1}}{I_{n}}{I_{n}}{J5}};

\node at ( $ (a)!0.5!(b) $ ) {$\longmapsto$};
\node at ( $ (b)!0.5!(c) $ ) {$\longmapsto$};
\node at ( $ (c)!0.5!(d) $ ) {\rotatebox{-90}{$\longmapsto$}};
\node at ( $ (d)!0.47!(e) $ ) {\rotatebox{180}{$\longmapsto$}};
\end{tikzpicture}
\]
This shows that the square matrix $S$ does transform by \emph{similarity}, i.e., by using transformations of the form
\[
S \longmapsto XSX^{-1},
\]
for non-singular matrices $X$ of the same order as $S$. Since the indecomposability of $M_{V}$ implies that of $S$, we obtain precisely the indecomposable representations of type $\mathrm{0}$.

\item At least one of the operators $f_{\alpha}$, $f_{\beta}$, $f_{\gamma}$, and $f_{\delta}$ is not an isomorphism. Since we are working up to duality, we can assume, without loss of generality, that 
$\Ima f_{\gamma}\neq V_{1}$. We proceed by induction on $n=\dim_{k}  V_{1}$. The case $n=1$ produces the matrix presentations 
\[
\renewcommand\arraystretch{1.45}
\begin{array}{|
  >{\centering\arraybackslash}p{\widthof{$\ile[0]$}}|
  >{\centering\arraybackslash}p{\widthof{$\ile[0]$}}|
  }
\hline
$1$ & $0$ \\
\hline
$1$ & $1$ \\
\hline
\end{array}\,,
\qquad
\begin{array}{|
  >{\centering\arraybackslash}p{\widthof{$\ile[0]$}}|
  >{\centering\arraybackslash}p{\widthof{$\ile[0]$}}|
  }
\hline
$1$ & $\ido[0]$ \\
\hline
$\ile[0]$ & $I_{0}$ \\
\hline
\end{array}\,,%
\qquad\text{and}\qquad
\begin{array}{|
  >{\centering\arraybackslash}p{\widthof{$\ile[0]$}}|
  >{\centering\arraybackslash}p{\widthof{$\ile[0]$}}|
  }
\hline
$1$ & $\iup[0]$ \\
\hline
$1$ & $\ido[0]$ \\
\hline
\end{array}
\]
which are of types $\mathrm{I}$, $\mathrm{II}$, and $\mathrm{IV}$, respectively.

Let $n>1$. We begin by placing at the bottom of $M_{12}$ the non-zero rows corresponding to $\Ima f_{\gamma}$:

\[
\begin{tikzpicture}[baseline=(current bounding box.center)]
\matrix[matrix of math nodes,
    ampersand replacement=\&,
    text width=3em,
    text height=1.7em,
    text depth=1em,
    align=center
]
(mim)
{
|[text height=0.7em,text depth=0.2em]|M'_{11} 
  \& |[text height=0.7em,text depth=0.2em]|0 \\
|[text height=1em,text depth=0.4em]|M''_{11} 
  \& |[text height=1em,text depth=0.4em]|M'_{12} 
\\
M_{21} \& M_{22} \\
};
\draw (mim-1-1.north west) -- (mim-3-1.south west);
\draw (mim-1-1.north east) -- (mim-3-1.south east);
\draw (mim-1-2.north east) -- (mim-3-2.south east);

\draw (mim-1-1.north west) -- (mim-1-2.north east);
\draw (mim-2-1.south west) -- (mim-2-2.south east);
\draw (mim-3-1.south west) -- (mim-3-2.south east);

\draw[dashed] (mim-1-1.south west) -- (mim-1-2.south east);

\draw[decoration={brace,mirror,raise=5pt},decorate]
  (mim-1-1.north west) --
    node[xshift=-7pt,anchor=east] {\footnotesize$V_{1}$}
  (mim-2-1.south west);
\draw[decoration={brace,raise=5pt},decorate]
  (mim-2-2.north east) --
    node[xshift=7pt,anchor=west] {\footnotesize$\Ima f_{\gamma}$}
  (mim-2-2.south east);

\end{tikzpicture}.
\] 
For the block $M'_{11}$ we can perform arbitrary elementary row and column operations, so $M'_{11}$ can be transformed into its standard form $S(M'_{11},r)$, where $r=\rank M'_{11}\geq 1$. The lower horizontal stripe of zeros in $S(M'_{11},r)$ has to be void; otherwise,  one could split off null direct summands from $V$, which contradicts its indecomposability. Therefore, using the identity $I_{r}$ to annihilate all entries  below it and inside the stripe $M''_{11}$, we obtain that $M_{V}$ is equivalent to a matrix of the form
\begin{equation}
\label{equ:redmat}
\begin{tikzpicture}[baseline=(current bounding box.center)]
\matrix[matrix of math nodes,
    ampersand replacement=\&,
    text width=2.5em,
    text height=1.7em,
    text depth=1em,
    align=center
]
(mim)
{
|[text width=0.8em,text height=0.7em,text depth=0em]|I_{r} 
  \&
|[text height=0.7em,text depth=0em,text width=1.8em]|0 
  \& |[text height=0.7em,text depth=0em]|0 
\\
|[text width=0.8em,text height=1em,text depth=0.4em]|0 
  \&
  |[text height=1em,text depth=0.4em,text width=1.8em]|M'''_{11} 
  \& |[text height=1em,text depth=0.4em]|M'_{12} 
\\
|[text width=0.8em]|L 
  \& |[text width=1.8em]|M'_{21} \& M_{22} \\
};
\draw (mim-1-1.north west) -- (mim-3-1.south west);
\draw (mim-1-2.north east) -- (mim-3-2.south east);
\draw (mim-1-3.north east) -- (mim-3-3.south east);
\draw[dashed] (mim-1-1.north east) -- (mim-3-1.south east);

\draw (mim-1-1.north west) -- (mim-1-3.north east);
\draw (mim-2-1.south west) -- (mim-2-3.south east);
\draw (mim-3-1.south west) -- (mim-3-3.south east);
\draw[dashed] (mim-1-1.south west) -- (mim-1-3.south east);

\node[anchor=west,outer sep=0pt,inner sep=0pt] at (mim.east) {\small, $r\geq 1$.};
\end{tikzpicture}
\end{equation}
Now, in the block $M'_{21}$ we can perform arbitrary elementary row and column operations, so this block can be reduced to its standard form $S(M'_{21},m)$, with $m\geq 1$. The identity $I_{m}$ can be then used to annihilate, by using suitable elementary column operations, the first $m$ rows of the block $L$ to its left. This leaves us with a matrix of the form 
\begin{equation}
\label{equ:redmat2}
\begin{tikzpicture}[baseline=(current bounding box.center)]
\matrix[matrix of math nodes,
    ampersand replacement=\&,
    text width=2.5em,
    text height=1.7em,
    text depth=1em,
    align=center
]
(mim)
{
|[text width=1em,text height=0.7em,text depth=0em]|I_{r} 
  \&
  |[text width=1.2em,text height=0.7em,text depth=0em]|0
  \& 
  |[text height=0.7em,text depth=0em,text width=1em]|0 
  \& |[text height=0.7em,text depth=0em]|0 
\\
|[text width=1em,text height=1em,text depth=0.4em]|0 
  \&
  |[text width=1.2em,text height=0.7em,text depth=0em]|
  \&  
  |[text height=1em,text depth=0.4em,text width=1em]| 
  \& |[text height=1em,text depth=0.4em]|M'_{12} 
\\
|[text width=1em,text height=1em,text depth=0.4em]|0 
  \&
  |[text width=1.2em,text height=0.7em,text depth=0em]|I_{m}
  \&  
  |[text height=1em,text depth=0.4em,text width=1em]|0 
  \& |[text height=1em,text depth=0.4em]|
\\
|[text width=1em,text height=1em,text depth=0.4em]|L' 
  \&
  |[text width=1.2em,text height=0.7em,text depth=0em]|0
  \&  
  |[text height=1em,text depth=0.4em,text width=1em]|0 
  \& |[text height=1em,text depth=0.4em]| 
\\
};
\draw (mim-1-1.north west) -- (mim-4-1.south west);
\draw (mim-1-3.north east) -- (mim-4-3.south east);
\draw (mim-1-4.north east) -- (mim-4-4.south east);

\draw[dashed] (mim-1-1.north east) -- (mim-4-1.south east);
\draw[dashed] (mim-1-3.north west) -- (mim-4-3.south west);
\draw (mim-1-1.north west) -- (mim-1-4.north east);
\draw (mim-2-1.south west) -- (mim-2-4.south east);
\draw (mim-4-1.south west) -- (mim-4-4.south east);
\draw[dashed] (mim-1-1.south west) -- (mim-1-4.south east);
\draw[dashed] (mim-3-1.south west) -- (mim-3-3.south east);
\node[fill=white] at (mim-2-2.east|-mim-2-4.center) {$M''''_{11}$};
\node[fill=white] at ([yshift=-1ex]mim-2-4.center|-mim-3-2.south) {$M'_{22}$};
\node[anchor=west,outer sep=0pt,inner sep=0pt,align=center] at (mim.east) {,$\begin{aligned}m&\geq 1\\r&\geq 1\end{aligned}$};
\end{tikzpicture}
\end{equation}
Notice that, in this matrix, the vertical or the horizontal stripes of zeros from $S(M'_{21},m)$ can be void. Finally, for the block $L'$ we can perform arbitrary elementary row and column operations. The latter will modify the block $I_{r}$ at the upper left corner in matrix~\eqref{equ:redmat2}, but the identity can be restored by applying the inverse operations on rows. This implies that $L'$ can be transformed into its standard form $S(L',s)$, for $r\geq s\geq 1$. Necessarily, we must have $r=s$ and the vertical stripe of zeros in $S(L',s)$ has to be void. Otherwise, direct summands of the form  $\begin{array}{|c|c|}\hline 1 & 0 \\ \hline 0 & 0 \\ \hline\end{array}$ would split off from the matrix~\eqref{equ:redmat2}, and this would contradict the indecomposability of $V$; this means that $M_{V}$ is equivalent to a matrix of the form 
\begin{equation}
\label{equ:redmat3}
\begin{tikzpicture}[baseline=(current bounding box.center)]
\matrix[matrix of math nodes,
    ampersand replacement=\&,
    text width=2.5em,
    text height=1.7em,
    text depth=1em,
    align=center
]
(mim)
{
|[text width=0.8em,text height=0.7em,text depth=0em]|I_{r} 
  \&
  |[text width=1.2em,text height=0.7em,text depth=0em]|0
  \& 
  |[text height=0.7em,text depth=0em,text width=1em]|0 
  \& |[text height=0.7em,text depth=0em]|0 
\\
|[text width=0.8em,text height=1em,text depth=0.4em]|0 
  \&
  |[text width=1.2em,text height=0.7em,text depth=0em]|
  \&  
  |[text height=1em,text depth=0.4em,text width=1em]| 
  \& |[text height=1em,text depth=0.4em]|M'_{12} 
\\
|[text width=0.8em,text height=0.8em,text depth=0.1em]|0 
  \&
  |[text width=1.2em,text height=0.8em,text depth=0.1em]|I_{m}
  \&  
  |[text height=0.8em,text depth=0.1em,text width=1em]|0 
  \& |[text height=0.8em,text depth=0.1em]|
\\
|[text width=0.8em,text height=0.8em,text depth=0.1em]|I_{r} 
  \&
  |[text width=1.2em,text height=0.8em,text depth=0.1em]|0
  \&  
  |[text height=0.8em,text depth=0.1em,text width=1em]|0 
  \& |[text height=0.8em,text depth=0.1em]| 
\\
|[text width=0.8em,text height=0.8em,text depth=0.1em]|0 
  \&
  |[text width=1em,text height=0.8em,text depth=0.1em]|0
  \&  
  |[text height=0.8em,text depth=0.1em,text width=1em]|0 
  \& |[text height=0.8em,text depth=0.1em]| 
\\
};
\draw (mim-1-1.north west) -- (mim-5-1.south west);
\draw (mim-1-3.north east) -- (mim-5-3.south east);
\draw (mim-1-4.north east) -- (mim-5-4.south east);

\draw[dashed] (mim-1-1.north east) -- (mim-5-1.south east);
\draw[dashed] (mim-1-3.north west) -- (mim-5-3.south west);
\draw (mim-1-1.north west) -- (mim-1-4.north east);
\draw (mim-2-1.south west) -- (mim-2-4.south east);
\draw (mim-5-1.south west) -- (mim-5-4.south east);
\draw[dashed] (mim-1-1.south west) -- (mim-1-4.south east);
\draw[dashed] (mim-3-1.south west) -- (mim-3-3.south east);
\draw[dashed] (mim-4-1.south west) -- (mim-4-3.south east);
\node[fill=white] at (mim-2-2.east|-mim-2-4.center) {$M'''''_{11}$};
\node[fill=white] at (mim-4-4.center) {$M''_{22}$};
\node[anchor=west,outer sep=0pt,inner sep=0pt,align=center] at (mim.east) {,$\begin{aligned}m&\geq 1\\r&\geq 1\end{aligned}$};
\end{tikzpicture}
\end{equation}
The last horizontal stripe of zeros and the rightmost vertical stripe of zeros of the lower left block can be void.
Notice now that the block matrices
\[
M_{V}=
\begin{tikzpicture}[baseline=(current bounding box.center)]
\matrix[matrix of math nodes,
    ampersand replacement=\&,
    text width=2em,
    text height=1.3em,
    text depth=0.6em,
    align=center
]
(mim)
{
M_{11} \& M_{12} \\
M_{21} \& M_{22} \\
};
\draw (mim-1-1.north west) -- (mim-2-1.south west);
\draw (mim-1-1.north east) -- (mim-2-1.south east);
\draw (mim-1-2.north east) -- (mim-2-2.south east);
\draw (mim-1-1.north west) -- (mim-1-2.north east);
\draw (mim-1-1.south west) -- (mim-1-2.south east);
\draw (mim-2-1.south west) -- (mim-2-2.south east);
\end{tikzpicture}
\qquad\text{and}\qquad
N=
\begin{tikzpicture}[baseline=(current bounding box.center)]
\matrix[matrix of math nodes,
    ampersand replacement=\&,
    text width=2.5em,
    text height=1.7em,
    text depth=1em,
    align=center
]
(mim)
{
|[text width=1.2em,text height=1.5em,text depth=0.7em]|
  \&  
  |[text height=1.5em,text depth=0.7em,text width=1em]| 
  \& |[text height=1.5em,text depth=0.7em]|M'_{12} 
\\
|[text width=1.2em,text height=0.8em,text depth=0.1em]|I_{m}
  \&  
  |[text height=0.8em,text depth=0.1em,text width=1em]|0 
  \& |[text height=0.8em,text depth=0.1em]|
\\
|[text width=1.2em,text height=0.8em,text depth=0.1em]|0
  \&  
  |[text height=0.8em,text depth=0.1em,text width=1em]|0 
  \& |[text height=0.8em,text depth=0.1em]| 
\\
|[text width=1.2em,text height=0.8em,text depth=0.1em]|0
  \&  
  |[text height=0.8em,text depth=0.1em,text width=1em]|0 
  \& |[text height=0.8em,text depth=0.1em]| 
\\
};
\draw (mim-1-1.north west) -- (mim-4-1.south west);
\draw (mim-1-2.north east) -- (mim-4-2.south east);
\draw (mim-1-3.north east) -- (mim-4-3.south east);
\draw[dashed] (mim-2-1.north east) -- (mim-4-1.south east);

\draw (mim-1-1.north west) -- (mim-1-3.north east);
\draw (mim-1-1.south west) -- (mim-1-3.south east);
\draw (mim-4-1.south west) -- (mim-4-3.south east);
\draw[dashed] (mim-2-1.south west) -- (mim-2-2.south east);
\draw[dashed] (mim-3-1.south west) -- (mim-3-2.south east);

\node at (mim-3-3.center) {$M''_{22}$};
\node at (mim-1-1.east|-mim-1-3.center) {$M'''''_{11}$};
\end{tikzpicture}
\]  
are matrix presentations of two representations which satisfy the conditions of Proposition~\ref{pro:reducedindecomposable}. We conclude that the matrix presentation $N$ is indecomposable and, thus, corresponds to an indecomposable representation 
\[
V'=(V'_{1},V'_{2},V'_{3},V'_{4},f'_{\alpha},f'_{\beta},f'_{\gamma},f'_{\delta})\]
with $\dim_{k} V'_{1}=n-r<n$. The induction hypothesis implies that its matrix presentation $N$ is equivalent, up to permutations of the vertical or horizontal stripes, to one of the forms presented in Figure~\ref{fig:indecomposableA3}. In order to finish our proof, all that remains to be done is to replace, in the matrix~\eqref{equ:redmat3}, the presentation corresponding to $N$ for each one of the forms from the list and then verify that their extended matrix forms~\eqref{equ:redmat3} again are in the list. The types $\mathrm{0}$, $\mathrm{III}^{\ast}$, and $\mathrm{IV}^{\ast}$ must be excluded from this process, since $\Codim_{k} V'_{i}=0$, for all $i\in\inte{4}$, leaving only the possibility for $V'$ to be of types $\mathrm{I}$, $\mathrm{II}$, $\mathrm{III}$ or $\mathrm{IV}$.  Therefore, all one has to do is to perform four direct restorations $V'\mapsto V$. Additionally, notice that for every matrix presentation in Figure~\ref{fig:indecomposableA3}, $\Codim_{k} V_{i}\leq 1$, for each $i\in\inte{4}$, i.e., $r=1$ and each of the identity blocks $I_{r}$ in the matrix~\eqref{equ:redmat3} is of order 1.

This means that the restoration process is really simple: in each of the four possible cases, it is necessary to take each of the matrix presentations $V'$ from Figure~\ref{fig:indecomposableA3} and, after permuting horizontal or vertical stripes, if required, place in the lower-left block a sub-block with codimension $1$. Then, add one new row to the top and one new column to the left. In the newly added column, two 1s are to be placed: one, in the topmost position (corresponding to the intersection of the newly added row and column) and  the other, in the new column, precisely at the position where a row  of zeros appears in the lower-right block of  $V'$.  All the other entries in these newly added rows and columns are zeros. Perform some final elementary operations with rows or columns and check that the resulting matrix also belongs to the list.\qedhere
\end{enumerate}
\end{proof}
\begin{rema}
We present here the exhaustive details of the reconstruction process $V'\mapsto V$ for the four possible types. In the diagrams that follow, an arrow of the form 
\tikz[>=Latex]{\coordinate (aux1); \coordinate[right=of aux1] (aux2);
\draw[<->]
  ([yshift=3pt]aux1)
    to[bend left]  
  ([yshift=3pt]aux2);
} indicates that the two vertical stripes must be swapped, while an arrow of the form
\tikz[>=Latex,baseline=(current bounding box.center)]{\coordinate (aux1); \coordinate[below=of aux1] (aux2);
\draw[<->]
  ([yshift=3pt]aux1)
    to[bend left]  
  ([yshift=3pt]aux2);
} signifies that the two horizontal stripes must be interchanged. An arrow of the form 
\tikz[>=Latex]{\coordinate (aux1); \coordinate[right=of aux1] (aux2);
\draw[->]
  ([yshift=3pt]aux1)
    to[bend left]  
  ([yshift=3pt]aux2);
} indicates that the column where it originates must be moved to the far right position within the same vertical stripe, using appropriate column swaps. Similarly, an arrow of the form
\tikz[>=Latex,baseline=(current bounding box.center)]{\coordinate (aux1); \coordinate[below=of aux1] (aux2);
\draw[->]
  ([yshift=3pt]aux1)
    to[bend left]  
  ([yshift=3pt]aux2);
} means that the row from which it originates has to be moved to the bottom position inside the same horizontal stripe, performing appropriate row swaps.

\begin{enumerate}[leftmargin=*]
\item For a type $\mathrm{I}$ matrix:\nopagebreak
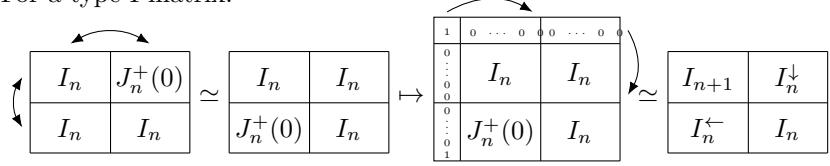
\begin{center}
\recmatrix[text width=3em,inner sep=0pt]{I_{n}}{J^{+}_{n}(0)}{I_{n}}{I_{n}}{matI1}
$\simeq$
\recmatrix[text width=3em,inner sep=0pt]{I_{n}}{I_{n}}{J^{+}_{n}(0)}{I_{n}}{matI2}
$\mapsto$\hspace*{10pt}
\recmatrix[text width=3em,inner sep=0pt,text depth=0.8em,text height=1.4em]{I_{n}}{I_{n}}{J^{+}_{n}(0)}{I_{n}}{matI3}
\begin{tikzpicture}[remember picture,overlay,>=Latex]
\draw 
  (matI3.south west) -- 
  ([xshift=-10pt]matI3.south west) -- 
  ([xshift=-10pt,yshift=10pt]matI3.north west) --
  ([yshift=10pt]matI3.north east) --
  (matI3.north east);
\draw 
  (matI3-1-1.south west) -- ++(-10pt,0pt);  
\draw 
  (matI3-1-1.north west) -- ++(-10pt,0pt);  
\draw 
  (matI3-1-1.north west) -- ++(0pt,10pt);  
\draw 
  (matI3-1-1.north east) -- ++(0pt,10pt);
\node[font=\tiny] 
  at ([xshift=-5pt,yshift=4.5pt]matI3-1-1.north west) 
  {\scalebox{.75}{$1$}};    
\node[anchor=north,font=\tiny] 
  at ([xshift=-5pt,yshift=2pt]matI3-1-1.north west) 
  {\scalebox{.75}{$\begin{matrix}0\\[-2ex]\vdotss\\0\\0\end{matrix}$}};    
\node[anchor=north,font=\tiny] 
  at ([xshift=-5pt,yshift=2pt]matI3-2-1.north west) 
  {\scalebox{.75}{$\begin{matrix}0\\[-2ex]\vdotss\\0\\1\end{matrix}$}};    
\node[anchor=west,font=\tiny] 
  at ([xshift=0pt,yshift=4pt]matI3-1-1.north west) 
  {\scalebox{.75}{$\renewcommand\arraycolsep{3pt}\begin{matrix}0 & \cdots & 0 & 0\end{matrix}$}};    
\node[anchor=west,font=\tiny] 
  at ([xshift=0pt,yshift=4pt]matI3-1-2.north west) 
  {\scalebox{.75}{$\renewcommand\arraycolsep{3pt}\begin{matrix}0 & \cdots & 0 & 0\end{matrix}$}};    
\draw[->] 
  ([xshift=-5pt,yshift=12pt]matI3-1-1.north west) 
    to[bend left]
  ([xshift=-3pt,yshift=12pt]matI3-1-1.north east);
\draw[->] 
  ([xshift=3pt,yshift=5pt]matI3-1-2.north east) 
    to[bend left]
  ([xshift=3pt,yshift=3pt]matI3-1-2.south east);
  
\draw[<->]
  ([yshift=3pt]matI1-1-1.north)
    to[bend left]  
  ([yshift=3pt]matI1-1-2.north);
\draw[<->]
  ([xshift=-3pt]matI1-1-1.west)
    to[bend right]  
  ([xshift=-3pt]matI1-2-1.west);
\end{tikzpicture}
$\simeq$
\recmatrix[text width=3em,inner sep=0pt]{I_{n+1}}{\ido}{\ile}{I_{n}}{matI4}
\end{center}

\item For a type $\mathrm{II}$ matrix:
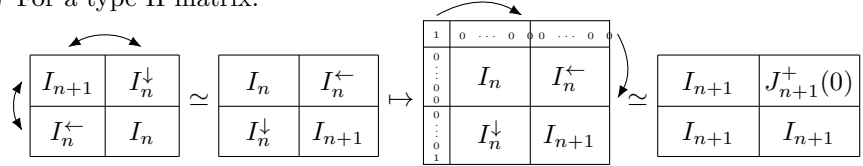
\begin{center}
\recmatrix[text width=2.8em,inner sep=0pt]{I_{n+1}}{\ido}{\ile}{I_{n}}{matII1}
$\simeq$
\recmatrix[text width=3em,inner sep=0pt]{I_{n}}{\ile}{\ido}{I_{n+1}}{matII2}
$\mapsto$\hspace*{10pt}
\recmatrix[text width=3em,inner sep=0pt,text depth=0.8em,text height=1.4em]{I_{n}}{\ile}{\ido}{I_{n+1}}{matII3}
\begin{tikzpicture}[remember picture,overlay,>=Latex]
\draw 
  (matII3.south west) -- 
  ([xshift=-10pt]matII3.south west) -- 
  ([xshift=-10pt,yshift=10pt]matII3.north west) --
  ([yshift=10pt]matII3.north east) --
  (matII3.north east);
\draw 
  (matII3-1-1.south west) -- ++(-10pt,0pt);  
\draw 
  (matII3-1-1.north west) -- ++(-10pt,0pt);  
\draw 
  (matII3-1-1.north west) -- ++(0pt,10pt);  
\draw 
  (matII3-1-1.north east) -- ++(0pt,10pt);
\node[font=\tiny] 
  at ([xshift=-5pt,yshift=4.5pt]matII3-1-1.north west) 
  {\scalebox{.75}{$1$}};    
\node[anchor=north,font=\tiny] 
  at ([xshift=-5pt,yshift=2pt]matII3-1-1.north west) 
  {\scalebox{.75}{$\begin{matrix}0\\[-2ex]\vdotss\\0\\0\end{matrix}$}};    
\node[anchor=north,font=\tiny] 
  at ([xshift=-5pt,yshift=2pt]matII3-2-1.north west) 
  {\scalebox{.75}{$\begin{matrix}0\\[-2ex]\vdotss\\0\\1\end{matrix}$}};    
\node[anchor=west,font=\tiny] 
  at ([xshift=0pt,yshift=4pt]matII3-1-1.north west) 
  {\scalebox{.75}{$\renewcommand\arraycolsep{3pt}\begin{matrix}0 & \cdots & 0 & 0\end{matrix}$}};    
\node[anchor=west,font=\tiny] 
  at ([xshift=0pt,yshift=4pt]matII3-1-2.north west) 
  {\scalebox{.75}{$\renewcommand\arraycolsep{3pt}\begin{matrix}0 & \cdots & 0 & 0\end{matrix}$}};    
\draw[->] 
  ([xshift=-5pt,yshift=12pt]matII3-1-1.north west) 
    to[bend left]
  ([xshift=-3pt,yshift=12pt]matII3-1-1.north east);
\draw[->] 
  ([xshift=3pt,yshift=5pt]matII3-1-2.north east) 
    to[bend left]
  ([xshift=3pt,yshift=3pt]matII3-1-2.south east);
  
\draw[<->]
  ([yshift=3pt]matII1-1-1.north)
    to[bend left]  
  ([yshift=3pt]matII1-1-2.north);
\draw[<->]
  ([xshift=-3pt]matII1-1-1.west)
    to[bend right]  
  ([xshift=-3pt]matII1-2-1.west);
\end{tikzpicture}
$\simeq$
\recmatrix[text width=3.8em,inner sep=0pt]{I_{n+1}}{J^{+}_{n+1}(0)}{I_{n+1}}{I_{n+1}}{matII4}
\end{center}

\item For a type $\mathrm{III}$ matrix:
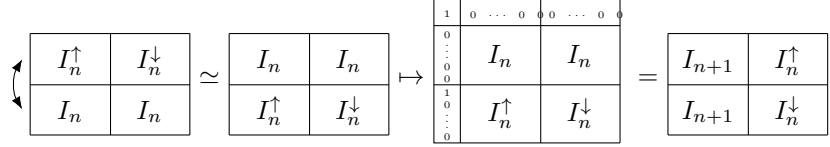
\begin{center}
\recmatrix[text width=3em,inner sep=0pt]{\iup}{\ido}{I_{n}}{I_{n}}{matIII1}
$\simeq$
\recmatrix[text width=3em,inner sep=0pt]{I_{n}}{I_{n}}{\iup}{\ido}{matIII2}
$\mapsto$\hspace*{10pt}
\recmatrix[text width=3em,inner sep=0pt,text depth=0.8em,text height=1.4em]{I_{n}}{I_{n}}{\iup}{\ido}{matIII3}
\begin{tikzpicture}[remember picture,overlay,>=Latex]
\draw 
  (matIII3.south west) -- 
  ([xshift=-10pt]matIII3.south west) -- 
  ([xshift=-10pt,yshift=10pt]matIII3.north west) --
  ([yshift=10pt]matIII3.north east) --
  (matIII3.north east);
\draw 
  (matIII3-1-1.south west) -- ++(-10pt,0pt);  
\draw 
  (matIII3-1-1.north west) -- ++(-10pt,0pt);  
\draw 
  (matIII3-1-1.north west) -- ++(0pt,10pt);  
\draw 
  (matIII3-1-1.north east) -- ++(0pt,10pt);
\node[font=\tiny] 
  at ([xshift=-5pt,yshift=4.5pt]matIII3-1-1.north west) 
  {\scalebox{.75}{$1$}};    
\node[anchor=north,font=\tiny] 
  at ([xshift=-5pt,yshift=2pt]matIII3-1-1.north west) 
  {\scalebox{.75}{$\begin{matrix}0\\[-2ex]\vdotss\\0\\0\end{matrix}$}};    
\node[anchor=north,font=\tiny] 
  at ([xshift=-5pt,yshift=2pt]matIII3-2-1.north west) 
  {\scalebox{.75}{$\begin{matrix}1\\0\\[-2ex]\vdotss\\0\end{matrix}$}};    
\node[anchor=west,font=\tiny] 
  at ([xshift=0pt,yshift=4pt]matIII3-1-1.north west) 
  {\scalebox{.75}{$\renewcommand\arraycolsep{3pt}\begin{matrix}0 & \cdots & 0 & 0\end{matrix}$}};    
\node[anchor=west,font=\tiny] 
  at ([xshift=0pt,yshift=4pt]matIII3-1-2.north west) 
  {\scalebox{.75}{$\renewcommand\arraycolsep{3pt}\begin{matrix}0 & \cdots & 0 & 0\end{matrix}$}};    
  
\draw[<->]
  ([xshift=-3pt]matIII1-1-1.west)
    to[bend right]  
  ([xshift=-3pt]matIII1-2-1.west);
\end{tikzpicture}
$=$
\recmatrix[text width=3em,inner sep=0pt]{I_{n+1}}{\iup}{I_{n+1}}{\ido}{matIII4}
\end{center}

\item For a type $\mathrm{IV}$ matrix:
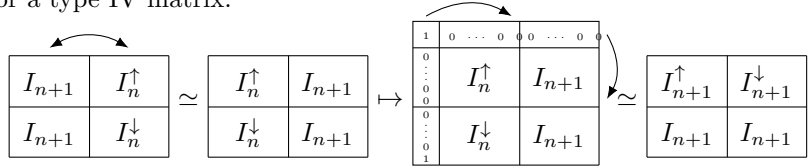
\begin{center}
\recmatrix[text width=3em,inner sep=0pt]{I_{n+1}}{\iup}{I_{n+1}}{\ido}{matIV1}
$\simeq$
\recmatrix[text width=3em,inner sep=0pt]{\iup}{I_{n+1}}{\ido}{I_{n+1}}{matIV2}
$\mapsto$\hspace*{10pt}
\recmatrix[text width=3em,inner sep=0pt,text depth=0.8em,text height=1.4em]{\iup}{I_{n+1}}{\ido}{I_{n+1}}{matIV3}
\begin{tikzpicture}[remember picture,overlay,>=Latex]
\draw 
  (matIV3.south west) -- 
  ([xshift=-10pt]matIV3.south west) -- 
  ([xshift=-10pt,yshift=10pt]matIV3.north west) --
  ([yshift=10pt]matIV3.north east) --
  (matIV3.north east);
\draw 
  (matIV3-1-1.south west) -- ++(-10pt,0pt);  
\draw 
  (matIV3-1-1.north west) -- ++(-10pt,0pt);  
\draw 
  (matIV3-1-1.north west) -- ++(0pt,10pt);  
\draw 
  (matIV3-1-1.north east) -- ++(0pt,10pt);
\node[font=\tiny] 
  at ([xshift=-5pt,yshift=4.5pt]matIV3-1-1.north west) 
  {\scalebox{.75}{$1$}};    
\node[anchor=north,font=\tiny] 
  at ([xshift=-5pt,yshift=2pt]matIV3-1-1.north west) 
  {\scalebox{.75}{$\begin{matrix}0\\[-2ex]\vdotss\\0\\0\end{matrix}$}};    
\node[anchor=north,font=\tiny] 
  at ([xshift=-5pt,yshift=2pt]matIV3-2-1.north west) 
  {\scalebox{.75}{$\begin{matrix}0\\[-2ex]\vdotss\\0\\1\end{matrix}$}};    
\node[anchor=west,font=\tiny] 
  at ([xshift=0pt,yshift=4pt]matIV3-1-1.north west) 
  {\scalebox{.75}{$\renewcommand\arraycolsep{3pt}\begin{matrix}0 & \cdots & 0 & 0\end{matrix}$}};    
\node[anchor=west,font=\tiny] 
  at ([xshift=0pt,yshift=4pt]matIV3-1-2.north west) 
  {\scalebox{.75}{$\renewcommand\arraycolsep{3pt}\begin{matrix}0 & \cdots & 0 & 0\end{matrix}$}};    
\draw[->] 
  ([xshift=-5pt,yshift=12pt]matIV3-1-1.north west) 
    to[bend left]
  ([xshift=-3pt,yshift=12pt]matIV3-1-1.north east);
\draw[->] 
  ([xshift=3pt,yshift=5pt]matIV3-1-2.north east) 
    to[bend left]
  ([xshift=3pt,yshift=3pt]matIV3-1-2.south east);
  
\draw[<->]
  ([yshift=3pt]matIV1-1-1.north)
    to[bend left]  
  ([yshift=3pt]matIV1-1-2.north);
\end{tikzpicture}
$\simeq$
\recmatrix[text width=3em,inner sep=0pt]{\iup[n+1]}{\ido[n+1]}{I_{n+1}}{I_{n+1}}{matIV4}
\end{center}
\end{enumerate}
For types $\mathrm{III}$ and $\mathrm{IV}$ there are two forms in which to perform the reconstruction process. In the diagrams above we presented one of them; the alternative ways will also produce the desired result. The details are left to the reader.
\end{rema}
\begin{rema}
The following diagram shows how the correspondence $V'\mapsto V$ transforms types $\mathrm{I}$, $\mathrm{II}$, $\mathrm{III}$, and $\mathrm{IV}$, as well as their corresponding dimension vectors:
\[
\begin{array}{c}
  \begin{tikzpicture}[node distance=1cm and 2cm]
  \node (I) {$\mathrm{I}$}; 
  \node[right=of I] (II) {$\mathrm{II}$}; 
  \draw[->] 
    (I.north east) to[bend left] 
    (II.north west);
  \draw[->] (II.south west) to[bend left] 
    (I.south east);
  \end{tikzpicture}
\\
  \begin{array}{ccc}
  (n,n,n,n)_{\mathrm{I}} & \mapsto & (n+1,n,n+1,n)_{\mathrm{II}} \\
  (n+1,n+1,n+1,n+1)_{\mathrm{I}} 
    & \rotatebox[origin=c]{180}{$\mapsto$} 
    &  \makebox[\widthof{$(n+1,n+1,n+1,n+1)_{\mathrm{I}}$}][c]{%
          $(n+1,n,n+1,n)_{\mathrm{II}}$%
          } \\
  \end{array}
\\[1.5em]
  \begin{tikzpicture}[node distance=1cm and 2cm]
  \node(III) {$\mathrm{III}$}; 
  \node[right=of III] (IV) {$\mathrm{IV}$}; 
  \draw[->] (III.north east) to[bend left] (IV.north west);
  \draw[->] (IV.south west) to[bend left] (III.south east);
  \end{tikzpicture}
\\
  \begin{array}{ccc}
  (n+1,n,n,n)_{\mathrm{III}} & \mapsto & (n+1,n+1,n+1,n)_{\mathrm{IV} } \\
  (n+2,n+1,n+1,n+1)_{\mathrm{III}} 
    & \rotatebox[origin=c]{180}{$\mapsto$} 
    & \makebox[\widthof{$(n+2,n+1,n+1,n+1)_{\mathrm{IV}}$}][c]{%
        $(n+1,n+1,n+1,n)_{\mathrm{IV}}$%
        }
  \end{array}
\end{array}
\]
The graphical invariants reflect these transformations. Indeed, to go from one type to another, as above, all one has to do is to keep one symbol in the upper right block fixed and move the other one to the block diagonally opposed:
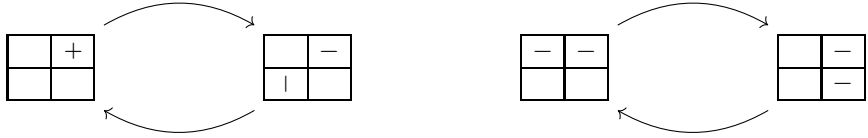
\[
  \begin{tikzpicture}[node distance=1cm and 2cm]
  \node (I) {\geomatrix{}{%
    \hspace{-3pt}{$-$}%
    \hspace*{\widthof{$-$}*\real{-0.5}}%
    \makebox[0pt][c]{\rotatebox[origin=c]{90}{$-$}}%
    }{}{}}; 
  \node[right=of I] (II) {\geomatrix{}{$-$}{\rotatebox[origin=c]{90}{$-$}}{}}; 
  \draw[->] 
    (I.north east) to[bend left] 
    (II.north west);
  \draw[->] (II.south west) to[bend left] 
    (I.south east);
  \node[right=2cm of II] (III) {\geomatrix{$-$}{$-$}{}{}};
  \node[right=of III] (IV) {\geomatrix{}{$-$}{}{$-$}};
  \draw[->] 
    (III.north east) to[bend left] 
    (IV.north west);
  \draw[->] (IV.south west) to[bend left] 
    (III.south east);
  \end{tikzpicture}
\]
\end{rema}
\begin{rema}
From the construction presented in the proof of Theorem~\ref{the:indecomposable2x2}, it follows that all the representations listed in Figure~\ref{fig:indecomposableA3} are indecomposable and mutually non-isomorphic. 
\end{rema}

\section{Endomorphism rings}
\label{sec:endomorphisms}
In this section we obtain our second main result: the classification of all the endomorphism rings for the indecomposable objects in $\repSk$. It is well known that given a representation $V$ of $\carc{S}$ over a field $k$, an \textit{endomorphism} of $V$ is a morphism $l\colon V\to V$ in $\repSk$. The set 
\[
\End V=\Hom(V,V)
=\{ l\colon V\to V \mid l\text{ is morphism}\}
\]
has a natural structure of ring with the usual sum and composition of morphisms. This is the \textit{ring of endomorphisms of V}. In the following theorem, we obtain these rings for all indecomposable representations of $\carc{S}$.
\begin{theo}
Let $V$ be an indecomposable representation with endomorphism ring $\mathcal{E}=\End V$. Then, the following holds:
\begin{enumerate}
\item\label{ite:endo0} If $V$ is of type $\mathrm{0}$, with $\dime=(n,n,n,n)$, then $\mathcal{E}\simeq k[t]/(p^{s}(t))$, where $p^{s}(t)$ is the minimal polynomial of the indecomposable Frobenius cell $F=F_{n}(p^{s}(t))$ (see Figure~\ref{fig:indecomposableA3}), with $p(t)\neq t$ and $s\deg p(t)=n$.

\item\label{ite:endo1} If $V$ is of type $\mathrm{I}$, with $\dime=(n,n,n,n)$, then $\mathcal{E}\simeq k[t]/(t^{n})$.

\item\label{ite:endo2} If $V$ is of type $\mathrm{II}$, with $\dime=(n+1,n,n+1,n)$, then $\mathcal{E}\simeq k[t]/(t^{n+1})$.

\item\label{ite:endo34} If $V$ is of type $\mathrm{III}$, $\mathrm{IV}$ or their duals, then $\mathcal{E}\simeq k$. 
\end{enumerate} 
\end{theo}
\begin{proof}
We will examine endomorphisms of an indecomposable representation $V=(V_{1},V_{1},V_{1},V_{1},f_{\alpha},f_{\beta},f_{\gamma},f_{\delta})$, that is, we will consider collections $l=(l_{1},l_{2},l_{3},l_{4})$ of linear maps  which make the following diagram commute:
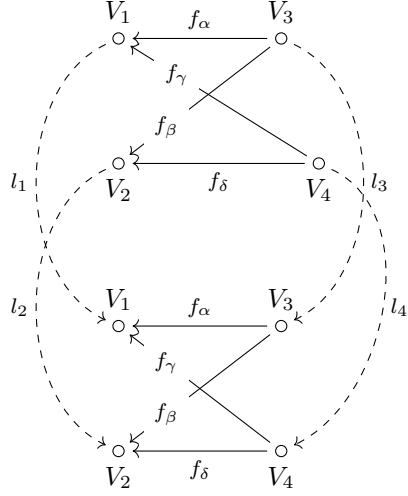
\begin{equation*}
\begin{tikzpicture}[
  baseline=(current bounding box.center),
  node distance=1.5cm and 2cm
  ]
\node[copo,label={above:$V_{1}$}] (v1) {};
\node[copo,below=of v1,label={below:$V_{2}$}] (v2) {};
\node[copo,right=of v1,label={above:$V_{3}$}] (v3) {};
\node[copo,below=of 3,label={below:$V_{4}$}] (v4) {};
\draw[->,shorten >= 3pt,shorten <= 3pt] 
  (v3) -- 
  node[above,font=\footnotesize] {$f_{\alpha}$} 
  (v1);
\draw[->,shorten >= 3pt,shorten <= 3pt] 
  (v3) -- 
  node[fill=white,pos=0.72] {\phantom{$f_{\beta}$}}
  node[font=\footnotesize,pos=0.72] {$f_{\beta}$}
  (v2);
\draw[->,shorten >= 3pt,shorten <= 3pt] 
  (v4) -- 
  node[fill=white,pos=0.72] {\phantom{$f_{\gamma}$}}
  node[font=\footnotesize,pos=0.72] {$f_{\gamma}$}
  (v1);
\draw[->,shorten >= 3pt,shorten <= 3pt] 
  (v4) -- 
  node[below,font=\footnotesize] {$f_{\delta}$}
  (v2);

\node[copo,below=2cm of v2,label={above:$V_{1}$}] (w1) {};
\node[copo,below=of w1,label={below:$V_{2}$}] (w2) {};
\node[copo,right=of w1,label={above:$V_{3}$}] (w3) {};
\node[copo,below=of w3,label={below:$V_{4}$}] (w4) {};
\draw[->,shorten >= 3pt,shorten <= 3pt] 
  (w3) -- 
  node[above,font=\footnotesize] {$f_{\alpha}$} 
  (w1);
\draw[->,shorten >= 3pt,shorten <= 3pt] 
  (w3) -- 
  node[fill=white,pos=0.72] {\phantom{$f_{\beta}$}}
  node[font=\footnotesize,pos=0.72] {$f_{\beta}$}
  (w2);
\draw[->,shorten >= 3pt,shorten <= 3pt] 
  (w4) -- 
  node[fill=white,pos=0.72] {\phantom{$f_{\gamma}$}}
  node[font=\footnotesize,pos=0.72] {$f_{\gamma}$}
  (w1);
\draw[->,shorten >= 3pt,shorten <= 3pt] 
  (w4) -- 
  node[below,font=\footnotesize] {$f_{\delta}$}
  (w2);

\foreach \index/\direct in {1,2}
{
\draw[->,shorten >= 3pt,shorten <= 3pt,dashed]
  (v\index) to[out=205,in=155]
  node[left,font=\footnotesize] {$l_{\index}$}
  (w\index);  
}
\foreach \index/\direct in {3,4}
{
\draw[->,shorten >= 3pt,shorten <= 3pt,dashed]
  (v\index) to[out=-25,in=25]
  node[right,font=\footnotesize] {$l_{\index}$}
  (w\index);  
}
\end{tikzpicture}
\end{equation*}
Equivalently, the following equalities must hold in order for $l$ to be an endomorphism of $V$:
\begin{equation}
\label{equ:endoequations}
l_{1}f_{\alpha}=f_{\alpha}l_{3},
\quad
l_{1}f_{\gamma}=f_{\gamma}l_{4},
\quad
l_{2}f_{\beta}=f_{\beta}l_{3},
\quad\text{and}\quad
l_{2}f_{\delta}=f_{\delta}l_{4}.
\end{equation}
We will work in terms of matrices. We fix ordered bases for spaces $V_{1},\ldots,V_{4}$ and denote by $[f_{\alpha}]$, $[f_{\beta}]$, $[f_{\gamma}]$, and $[f_{\delta}]$ the matrices of the corresponding operators with respect to the chosen bases. Similarly, $L_{1},\ldots,L_{4}$ will denote the matrices of the maps $l_{1},\ldots,l_{4}$, respectively.

\ref{ite:endo0} For a representation $V$ of type $\mathrm{0}$ with $\dime =(n,n,n,n)$, we have $V_{i}\simeq k^{n}$, for all $i\in\inte{4}$, $[f_{\alpha}]=[f_{\beta}]=[f_{\delta}]=I_{n}$, and $[f_{\gamma}]=F$. In this situation, the equations~\eqref{equ:endoequations} become:
\begin{equation*}
L_{1}I_{n}=I_{n}L_{3},
\quad
L_{1}F=FL_{4},
\quad
L_{2}I_{n}=I_{n}L_{3},
\quad\text{and}\quad
L_{2}I_{n}=I_{n}L_{4},
\end{equation*}
i.e., $L_{1}=L_{2}=L_{3}=L_{4}$ and $L_{1}F=FL_{1}$. Thus, in order to determine $\mathcal{E}$, we must identify the $k$-algebra 
\[
C(F)=\{L\in k^{n\times n} \mid LF=FL \},
\]
the center of $F$ in $k^{n\times n}$. Let $\alpha$ be a cyclic vector for $F$ and define the linear map $f\colon Z\to k^{n}$ by $f(L)=L\alpha$. We will prove that f is an isomorphism of \knob{vector} spaces. If $0=f(L)=L\alpha$, then $L(F^{i}\alpha)=F^{i}(L\alpha)=0$, for all $i\in\inte{n-1}$ and, hence, $L=0$ and $f$ is injective. Moreover, $\langle \{f(F^{j})\mid j\in\{0,\ldots,n-1\}\} \rangle=\langle \{F^{j}\alpha\mid j\in\{0,\ldots,n-1\}\}\rangle=k^{n}$, so $f$ is surjective. In particular, we have $\dim_{k}C(F)=n$. Since $\{I,F,\ldots,F^{n-1}\}$ is a basis for $k[F]$ as \knob{vector} space and clearly $k[F]\subseteq C(F)$, this implies that $C(F)=k[F]$. We have, thus, that the ring $C(F)$ is equal to
\[
\bigl\{ a_{0}F^{0} + a_{1}F^{1} + \cdots + a_{n-1}F^{n-1} 
\mid a_{i}\in k\text{, for all $i\in\{0,\ldots,n-1\}$} \bigr\}
\]
and, from here, $\mathcal{E}=k[t]/(p^{s}(t))$.

\ref{ite:endo1} For a representation $V$ of type $\mathrm{I}$ with $\dime =(n,n,n,n)$, we have the same situation as in~\ref{ite:endo0} for $F=J^{+}_{n}(0)$, so here $p^{s}(t)=t^{n}$.

\ref{ite:endo2} For a representation $V$ of type $\mathrm{II}$ with $\dime =(n+1,n,n+1,n)$, we have $V_{1}\simeq V_{3} \simeq k^{n+1}$, $V_{2}\simeq V_{4} \simeq k^{n}$, $[f_{\alpha}]=I_{n+1}$, $[f_{\beta}]=\ile$, $[f_{\gamma}]=\ido$, and $[f_{\delta}]=I_{n}$. Now, the equations~\eqref{equ:endoequations} become:
\begin{equation*}
L_{1}I_{n+1}=I_{n+1}L_{3},
\quad
L_{1}\ido=\ido L_{4},
\quad
L_{2}\ile=\ile L_{3},
\quad\text{and}\quad
L_{2}I_{n}=I_{n}L_{4},
\end{equation*}
i.e., $L_{1}=L_{3}$, $L_{2}=L_{4}$, 
\begin{equation}
\label{equ:type2}
L_{1}\ido = \ido L_{2}\quad\text{and}\quad
L_{2}\ile = \ile L_{1}.
\end{equation}
From~\eqref{equ:type2} we get
\begin{align*}
L_{1} &= a_{0}J^{0} + a_{1}J^{1} + \cdots + a_{n}J^{n}, \\
L_{2} &= a_{0}J^{0} + a_{1}J^{1} + \cdots + a_{n-1}J^{n-1},
\end{align*}
where $a_{i}\in k$, for all $i\in\{0,\ldots,n\}$ and $J=J^{+}_{n}(0)$. It is then clear that $\mathcal{E}\simeq k[t]/(t^{n+1})$.
 
\ref{ite:endo34} For a representation $V$ of type $\mathrm{III}$ with $\dime =(n+1,n,n,n)$, we have $V_{1} \simeq k^{n+1}$, $V_{2}\simeq V_{3}\simeq V_{4} \simeq k^{n}$, $[f_{\alpha}]=\iup$, $[f_{\beta}]=I_{n}$, $[f_{\gamma}]=\ido$, and $[f_{\delta}]=I_{n}$. In this case, the equations~\eqref{equ:endoequations} become:
\begin{equation*}
L_{1}\iup=\iup L_{3},
\quad
L_{1}\ido=\ido L_{4},
\quad
L_{2}I_{n}=I_{n} L_{3},
\quad\text{and}\quad
L_{2}I_{n}=I_{n}L_{4},
\end{equation*}
i.e., $L_{2}=L_{3}=L_{4}$ and
\begin{equation}
\label{equ:type3}
L_{1}\iup = \iup L_{2}\quad\text{and}\quad
L_{1}\ido = \ido L_{2}.
\end{equation}
From~\eqref{equ:type2} we get $L_{1}=aI_{n+1}$ and $L_{2}=aI_{n}$, for some $a\in k$. From this we conclude $\mathcal{E}\simeq k$.
 
For a representation $V$ of type $\mathrm{IV}$ with $\dime =(n+1,n+1,n+1,n)$, we have $V_{1} \simeq V_{2}\simeq V_{3}\simeq k^{n+1}$, $V_{4} \simeq k^{n}$, $[f_{\alpha}]=I_{n+1}$, $[f_{\beta}]=I_{n+1}$, $[f_{\gamma}]=\iup$, and $[f_{\delta}]=\ido$. Now, the equations~\eqref{equ:endoequations} become:
\begin{equation*}
L_{1}I_{n+1}=I_{n+1}L_{3},
\quad
L_{1}\iup=\iup L_{4},
\quad
L_{2}I_{n+1}=I_{n+1} L_{3},
\quad\text{and}\quad
L_{2}\ido =\ido L_{4},
\end{equation*}
i.e., $L_{1}=L_{2}=L_{3}$ and 
\begin{equation}
\label{equ:type4}
L_{1}\iup = \iup L_{4}\quad\text{and}\quad
L_{1}\ido = \ido L_{4}.
\end{equation}
From~\eqref{equ:type4} we get $L_{1}=aI_{n+1}$ and $L_{4}=aI_{n}$, for some $a\in k$. From this we conclude $\mathcal{E}\simeq k$.
 
Finally, for the representations of types $\mathrm{III}^{\ast}$ and $\mathrm{IV}^{\ast}$ all we have to do is to notice that one has the isomorphism
\[
\Hom(V^{\ast},V^{\ast})\simeq \Hom(V,V)
\]
and the result for the duals follows from what we just did for types $\mathrm{III}$ and $\mathrm{IV}$. This concludes our proof.
\end{proof}

\section{The quiver $\carc{S}$ and some other classification problems}
\label{sec:subproblems}
In this section we would like to illustrate how to use the list of matrix presentations given in Figure~\ref{fig:indecomposableA3}, to easily obtain solutions for similar classification problems. We will discuss the Kronecker problem and its contragredient variation.

The Kronecker problem is a well-known problem in the theory of linear transformations which involves the classification of all pairs of linear transformations between two finite-dimensional vector spaces over a field. Weierstrass provided a partial solution~\parencite{Wei1868}, while Kronecker offered a comprehensive solution~\parencite{Kro1990}. Other recent solutions, using a variety of techniques, can be found in, for example, \textcite{DorMed2023,Die46,Dev84,GabRoi1992,Zav2007}. 
A closely related variant of the Kronecker problem is its contragredient version, where the domain and codomain of one of the two operators are swapped. This problem was solved by Dobrovol'skaja and Ponomarev in \textcite{DobPon1965}. A more recent generalization was provided in \textcite{Zav2007}.
\begin{figure}[H]
  \centering
  \begin{subcaptionblock}{.48\textwidth}
  \centering
  \resizebox{.9\linewidth}{!}{%
  \begin{tabular}{@{}c@{\hspace{4.8em}}c@{}}
    \begin{tikzpicture}[
      node distance=1.5cm and 2cm,
      baseline=(current bounding box.center)]
    \node[copo,label={left:$V_{1}$}] (1) {};
    \node[copo,below=of 1,label={left:$V_{2}$}] (2) {};
    \node[copo,right=of 1,label={right:$V_{3}$}] (3) {};
    \node[copo,below=of 3,label={right:$V_{4}$}] (4) {};
    \draw[->,shorten >= 3pt,shorten <= 3pt] 
      (3) -- 
      node[above,font=\footnotesize] {$f_{\alpha}$} 
      (1);
    \draw[->,shorten >= 3pt,shorten <= 3pt] 
      (3) -- 
      node[fill=white,pos=0.7] {\phantom{$\beta\gamma$}}
      node[font=\footnotesize,pos=0.72] {$f_{\beta}=1$}
      (2);
    \draw[->,shorten >= 3pt,shorten <= 3pt] 
      (4) -- 
      node[fill=white,pos=0.7] {\phantom{$\gamma$}}
      node[font=\footnotesize,pos=0.7] {$f_{\gamma}$}
      (1);
    \draw[->,shorten >= 3pt,shorten <= 3pt] 
      (4) -- 
      node[below,font=\footnotesize] {$f_{\delta}=1$}
      (2);
    \end{tikzpicture}
  &
    \sqmatrix{A}{C}{B}{D}
  \\[8ex]
    \rotatebox[origin=c]{-90}{$\mapsto$} 
  &
    \rotatebox[origin=c]{-90}{$\mapsto$} 
  \\
    \begin{tikzpicture}[
      node distance=1.5cm and 2cm,
      baseline=(current bounding box.center)]
    \node[copo,label={left:$V_{1}$}] (1) {};
    \node[copo,right=of 1,label={right:$V_{3}$\rlap{${}\simeq V_{2}\simeq V_{4}$}}] (3) {};
    \draw[->,shorten >= 3pt,shorten <= 3pt] 
      ( $ (3) + (-3pt,3pt) $ ) -- 
      node[above,font=\footnotesize] {$f_{\alpha}$} 
      ( $ (1) + (3pt,3pt) $ );
    \draw[->,shorten >= 3pt,shorten <= 3pt] 
      ( $ (3) + (-3pt,-3pt) $ ) -- 
      node[below,font=\footnotesize] {$f_{\gamma}$}
      ( $ (1) + (3pt,-3pt) $ );
    \end{tikzpicture}
  &
    \sqmatrix{A}{C}{I}{I}  
  \end{tabular}}
  \caption{Starting from $\carc{S}$ and taking $f_{\beta}$ and $f_{\delta}$ as identities, we can identify $V_{4}$, $V_{3}$, and $V_{2}$ obtaining the Kronecker quiver}
  \label{sfig:A3tokronecker}
  \end{subcaptionblock}\hfill%
  \begin{subcaptionblock}{.48\textwidth}
  \centering
  \resizebox{.9\linewidth}{!}{%
  \begin{tabular}{@{}c@{\hspace{4em}}c@{}}
    \begin{tikzpicture}[
      node distance=1.5cm and 2cm,
      baseline=(current bounding box.center)]
    \node[copo,label={left:$V_{1}$}] (1) {};
    \node[copo,below=of 1,label={left:$V_{2}$}] (2) {};
    \node[copo,right=of 1,label={right:$V_{3}$}] (3) {};
    \node[copo,below=of 3,label={right:$V_{4}$}] (4) {};
    \draw[->,shorten >= 3pt,shorten <= 3pt] 
      (3) -- 
      node[above,font=\footnotesize] {$f_{\alpha}$} 
      (1);
    \draw[->,shorten >= 3pt,shorten <= 3pt] 
      (3) -- 
      node[fill=white,pos=0.7] {\phantom{$\beta\gamma$}}
      node[font=\footnotesize,pos=0.72] {$f_{\beta}=1$}
      (2);
    \draw[->,shorten >= 3pt,shorten <= 3pt] 
      (4) -- 
      node[fill=white,pos=0.7] {\phantom{$\beta\gamma$}}
      node[font=\footnotesize,pos=0.72] {$f_{\gamma}=1$}
      (1);
    \draw[->,shorten >= 3pt,shorten <= 3pt] 
      (4) -- 
      node[below,font=\footnotesize] {$f_{\delta}$}
      (2);
    \end{tikzpicture}
  &
    \sqmatrix{A}{C}{B}{D}
  \\[8ex]
    \rotatebox[origin=c]{-90}{$\mapsto$} 
  &
    \rotatebox[origin=c]{-90}{$\mapsto$} 
  \\
    \begin{tikzpicture}[
      node distance=1.5cm and 2cm,
      baseline=(current bounding box.center)]
    \node[copo,label={left:\llap{$V_{4}\simeq{}$}$V_{1}$}] (1) {};
    \node[copo,right=of 1,label={right:$V_{3}$\rlap{${}\simeq V_{2}$}}] (3) {};
    \draw[->,shorten >= 3pt,shorten <= 3pt] 
      ( $ (3) + (-3pt,3pt) $ ) -- 
      node[above,font=\footnotesize] {$f_{\alpha}$} 
      ( $ (1) + (3pt,3pt) $ );
    \draw[<-,shorten >= 3pt,shorten <= 3pt] 
      ( $ (3) + (-3pt,-3pt) $ ) -- 
      node[below,font=\footnotesize] {$f_{\delta}$}
      ( $ (1) + (3pt,-3pt) $ );
    \end{tikzpicture}
  &
    \sqmatrix{A}{I}{I}{D}  
  \end{tabular}}
  \caption{Starting from $\carc{S}$ and taking $f_{\beta}$ and $f_{\gamma}$ as identities, we can identify $V_{4}$ with $V_{1}$, and $V_{3}$ with $V_{2}$, obtaining the contragredient Kronecker quiver}
  \label{sfig:A3tocontragredient}
  \end{subcaptionblock}%
  \caption{Depiction of the identifications and matrix transformations that show that the classification problem for $\carc{S}$ contains as sub-problems the Kronecker problem and its contragredient version.}
  \label{fig:A3subproblems}
\end{figure}
In Figure~\ref{fig:A3subproblems}, we present some graphical schemes that show that the classical Kronecker problem and its contragredient version are embedded as sub-problems into the classification problem of the quiver $\carc{S}$. In fact, the matrix transformations that accompany each diagram also show how to obtain the solution to each subproblem from the solution given for $\carc{S}$ in Figure~\ref{fig:indecomposableA3}. We encourage the reader to fill up the details of the corresponding solutions.
\begin{coro}
\label{cor:kronecker}
All the indecomposable matrix presentations of the Kronecker problem are exhausted, up to isomorphism, by the presentations of the four types shown below in matrix form:\par\smallskip
\begin{tabular}{@{\hspace{1.5em}}>{$}r<{$} >{$}l<{$}}
\mathrm{0}=\mathrm{0}^{\ast}:
&  \hormatrix[
     baseline={([yshift=-0.75ex]current bounding box.center)}
     ]
     {I_{n}}{F}
     \text{ with
     $F=F_{n}(p^{s}(t)), p(t)\neq t$, and $n\geq 1$.} 
\\[2.5ex]
\mathrm{I}=\mathrm{I}^{\ast}:
&  \hormatrix[
     remember picture,
     baseline={([yshift=-0.75ex]current bounding box.center)}
     ]
    {I_{n}}{J^{+}_{n}(0)}%
    \text{ and }%
    \hormatrix[
    remember picture,
     baseline={([yshift=-0.75ex]current bounding box.center)}
    ]
    {J^{+}_{n}(0)}{I_{n}} \text{\,, for $n\geq 1$.}
\\[2.5ex]
\mathrm{II}=\mathrm{III}^{\ast}:
&  \hormatrix[
     remember picture,
     baseline={([yshift=-0.75ex]current bounding box.center)}
     ]
     {\ido}{\iup} \text{\,, for $n\geq 0$.}
\\[2.5ex]
\mathrm{III}=\mathrm{II}^{\ast}:
&  \hormatrix[
     remember picture,
     baseline={([yshift=-0.75ex]current bounding box.center)}
     ]
     {\iri}{\ile} \text{\,, for $n\geq 0$.}
\end{tabular}\par\medskip\noindent
The given types are pair-wise non-equivalent.\hfill\qedsymbol
\end{coro}
\begin{coro}
\label{cor:contragredient}
All the indecomposable matrix presentations of the contragredient Kronecker problem are exhausted, up to isomorphism, by the presentations of the four types shown below in matrix form:\par\smallskip
\begin{tabular}{@{\hspace{1.5em}}>{$}r<{$} >{$}l<{$}}
\mathrm{0}=\mathrm{0}^{\ast}:
&  \hormatrix[
     baseline={([yshift=-0.75ex]current bounding box.center)}
     ]
     {I_{n}}{F}
     \text{ with
     $F=F_{n}(p^{s}(t)), p(t)\neq t$, and $n\geq 1$.} 
\\[2.5ex]
\mathrm{I}=\mathrm{I}^{\ast}:
&  \hormatrix[
     remember picture,
     baseline={([yshift=-0.75ex]current bounding box.center)}
     ]
    {I_{n}}{J^{+}_{n}(0)}%
    \text{ and }%
    \hormatrix[
    remember picture,
     baseline={([yshift=-0.75ex]current bounding box.center)}
    ]
    {J^{+}_{n}(0)}{I_{n}} \text{\,, for $n\geq 1$.}
\\[2.5ex]
\mathrm{II}=\mathrm{III}^{\ast}:
&  \hormatrix[
     remember picture,
     baseline={([yshift=-0.75ex]current bounding box.center)}
     ]
     {\iup}{\iri} \text{\,, for $n\geq 0$.}
\\[2.5ex]
\mathrm{III}=\mathrm{II}^{\ast}:
&  \hormatrix[
     remember picture,
     baseline={([yshift=-0.75ex]current bounding box.center)}
     ]
     {\iri}{\iup} \text{\,, for $n\geq 0$.}
\end{tabular}\par\medskip\noindent
The given types are pair-wise non-equivalent.\hfill\qedsymbol
\end{coro}
\begin{rema}
As the reader will have noticed, the Kronecker problem and its contragredient version have very similar statements and solutions. This fact can be formalized by defining some appropriate functors between the categories involved and the category $\repSk$.

The classification problem for the quiver $\carc{S}$ contains three other sub-problems: the problem associated to a quiver of type $\widetilde{A}_{2}$; the classification problem for a linear relation, and the classification problem for pairs of linear relations. However, those sub-problems are beyond the scope of this paper and will be addressed elsewhere.
\end{rema}

\appendix
\section*{Appendix}
Each matrix presentation of the types $\mathrm{I}$, $\mathrm{II}$, $\mathrm{III}$, $\mathrm{IV}$, $\mathrm{III}^{\ast}$, and $\mathrm{IV}^{\ast}$ contains at most two 1s in each row and column. If a row (or column) contains two 1s, each of them appears in different adjacent blocks.  All other entries in the matrix are zeros. In the corresponding graphical invariants, the symbol $-$ (or $\rotatebox[origin=c]{90}{$-$}$) indicates a horizontal (or vertical) stripe containing a row (or column) with a single 1. This row (or column) has zeros in the block containing the symbol, and the 1 is located in the horizontal (or vertical) adjacent block. There are three mutually exclusive possibilities:
\begin{enumerate}[nolistsep]
\item The matrix has exactly one row and exactly one column, each containing a single 1. This situation corresponds to the types $\mathrm{I}$ and $\mathrm{II}$.
\item The matrix has exactly two rows, each containing a single 1. This situation corresponds to the types $\mathrm{III}$ and $\mathrm{IV}$.
\item The matrix has exactly two columns, each containing a single 1. This situation corresponds to the types $\mathrm{III}^{\ast}$ and $\mathrm{IV}^{\ast}$.
\end{enumerate}
All the 1s in a matrix presentation of the mentioned types can be \textquote{linked} together along a spiral-like path that starts at one of the special symbols and ends at the other. The algorithm that we present below implements this process. 

\begin{algorithm}[H]
\caption{Producing the matrix presentations for indecomposable objects from their graphical invariants}
\label{alg:graphicinvariant}
\begin{algorithmic}[1]
\Require{$type\in\{1,2,3,3^{\ast},4,4^{\ast}\}$ and $n$ is an integer $n\geq 1$}
\Statex
\If{$type = 1$}
  \Let{$num\_ones$}{$4n-1$}\Comment{$num\_ones$ is the total number of $1$s for a given $type$ and $n$}
\ElsIf{$type = 2$}
  \Let{$num\_ones$}{$4n+1$}
\ElsIf{$type = 3$ or $type = 4^{\ast}$}
  \Let{$num\_ones$}{$4n$}
\ElsIf{$type = 4$ or $type = 3^{\ast}$}
  \Let{$num\_ones$}{$4n+2$}
\EndIf
\Statex
\Let{$vxpos$}{1}\Comment{$vxpos$ is a valid $x$ coordinate to place a $1$}
\Let{$vypos$}{1}\Comment{$vypos$ is a valid $y$ coordinate to place a $1$}
\If{$type = 1$ or $type = 2$ or $type = 3$ or $type = 4$}
  \Let{$factorx\_1$}{$-1$}
  \Let{$factory\_1$}{$1$}
  \Let{$factorx\_3$}{$1$}
  \Let{$factory\_3$}{$-1$}
  \Let{$ff\_1$}{$2$}
  \Let{$fs\_2$}{$4$}
  \Let{$ff\_3$}{$4$}
  \Let{$fs\_4$}{$2$}
  \Forr{i\gets 1}{n}
    \Let{$m[i][1]$}{$0$}
  \EndFor
\ElsIf{$type = 3^{\ast}$ or $type = 4^{\ast}$}
  \Let{$factorx\_1$}{$1$}
  \Let{$factory\_1$}{$-1$}
  \Let{$factorx\_3$}{$-1$}
  \Let{$factory\_3$}{$1$}
\algstore{myalgo}
\end{algorithmic}
\end{algorithm}
\begin{algorithm}
\begin{algorithmic}[1]
\algrestore{myalgo}
  \Let{$ff\_1$}{$4$}
  \Let{$fs\_2$}{$2$}
  \Let{$ff\_3$}{$2$}
  \Let{$fs\_4$}{$4$}
  \Forr{j\gets 1}{n}
    \Let{$m[1][j]$}{$0$}
  \EndFor
\EndIf
\Statex
\Forr{i\gets 1}{num\_ones}\Comment{In each iteration, a 1 is placed and then the intermediate positions are filled with zeros}
  \If{$i\equiv 1\bmod 4$}
    \Let{$m[factorx\_1*vxpos][factory\_1*vypos]$}{$1$}
    \State $\Call{fill}{ff\_1,3}$  
  \ElsIf{$i\equiv 2\bmod 4$}
    \Let{$m[-vxpos][-vypos]$}{$1$}
    \State $\Call{fill}{3,fs\_2}$  
  \ElsIf{$i\equiv 3\bmod 4$}
    \Let{$m[factorx\_3*vxpos][factory\_3*vypos]$}{$1$}
    \State $\Call{fill}{ff\_3,1}$  
  \Else
    \Let{$m[vxpos][vypos]$}{$1$}
    \State $\Call{fill}{1,fs\_4}$  
  \EndIf
\EndFor

\State \Return $m$
\Ensure{$m$ is the matrix presentation of the indecomposable having the given type and size}
\algstore{myalgo}
\end{algorithmic}
\end{algorithm}

\begin{algorithm}
\begin{multicols}{2}
\begin{algorithmic}[1]
\algrestore{myalgo}
\Procedure{fill}{p,q}\Comment{Fills with zeros the appropriate positions in the quadrants $p$ and $q$}
\If{p=2 and q=3}
  \For{$k \gets 1 \textrm{ to } vypos-1$}
    \Let{$m[-vxpos][k]$}{$0$}
    \Let{$m[-vxpos][-k]$}{$0$}
  \EndFor
\ElsIf{p=3 and q=4}
  \For{$k \gets 1 \textrm{ to } vxpos-1$}
    \Let{$m[-k][-vypos]$}{$0$}
    \Let{$m[k][-vypos]$}{$0$}
  \EndFor
\ElsIf{p=4 and q=1}
  \For{$k \gets 1 \textrm{ to } vypos-1$}
    \Let{$m[vxpos][-k]$}{$0$}
  \EndFor
  \For{$k \gets 1 \textrm{ to } vypos$}
    \Let{$m[vxpos][k]$}{$0$}
  \EndFor
  \Let{$vypos$}{$vypos+1$}
\ElsIf{p=1 and q=2}
  \For{$k \gets 1 \textrm{ to } vxpos-1$}
    \Let{$m[k][vypos]$}{$0$}
  \EndFor
  \Let{$vxpos$}{$vxpos+1$}
  \For{$k \gets 1 \textrm{ to } vxpos-1$}
    \Let{$m[-k][vypos]$}{$0$}
  \EndFor
\ElsIf{p=4 and q=3}
  \For{$k \gets 1 \textrm{ to } vxpos-1$}
    \Let{$m[k][-vypos]$}{$0$}
    \Let{$m[-k][-vypos]$}{$0$}
  \EndFor
\ElsIf{p=3 and q=2}
  \For{$k \gets 1 \textrm{ to } vypos-1$}
    \Let{$m[-vxpos][-k]$}{$0$}
    \Let{$m[-vxpos][k]$}{$0$}
  \EndFor
\ElsIf{p=2 and q=1}
  \For{$k \gets 1 \textrm{ to } vxpos-1$}
    \Let{$m[-k][vypos]$}{$0$}
  \EndFor
  \For{$k \gets 1 \textrm{ to } vxpos$}
    \Let{$m[k][vypos]$}{$0$}
  \EndFor
  \Let{$vxpos$}{$vxpos+1$}
\ElsIf{p=1 and q=4}
  \For{$k \gets 1 \textrm{ to } vypos-1$}
    \Let{$m[vxpos][k]$}{$0$}
  \EndFor
  \Let{$vypos$}{$vypos+1$}
  \For{$k \gets 1 \textrm{ to } vypos-1$}
    \Let{$m[vxpos][-k]$}{$0$}
  \EndFor
\EndIf
\EndProcedure
\end{algorithmic}
\end{multicols}
\end{algorithm}

\FloatBarrier

The following diagram illustrates how the algorithm builds the matrix presentation in the case $type=2$ and $n=2$:
\[
\begin{tikzpicture}[baseline=(current bounding box.center)]
\matrix[mimat,text width=1.3em,text height=0.9em,text depth=0.2em] 
(mat)
{
1 \& 0 \& 0 \& 0 \& 1 \\
0 \& 1 \& 0 \& 1 \& 0 \\
0 \& 0 \& 1 \& 0 \& 0 \\
0 \& 0 \& 1 \& 1 \& 0 \\
0 \& 1 \& 0 \& 0 \& 1 \\
};
\draw (mat-1-1.north west) -- (mat-5-1.south west);
\draw (mat-1-4.north west) -- (mat-5-4.south west);
\draw (mat-1-5.north east) -- (mat-5-5.south east);

\draw (mat-1-1.north west) -- (mat-1-5.north east);
\draw (mat-4-1.north west) -- (mat-4-5.north east);
\draw (mat-5-1.south west) -- (mat-5-5.south east);

\begin{scope}[very thick,decoration={
    markings,
    }
    ] 
\tikzset{
myarrow/.style={
postaction={%
  decorate,
  decoration={markings,
  mark=at position #1 with {\arrow{>}}}}
  }
}  
\begin{pgfonlayer}{background}
\draw[opacity=0.4,myarrow=0.15] 
  ([xshift=25pt]mat-3-5.center) --
  (mat-3-3.center);
\draw[opacity=0.4,myarrow=0.4] 
  (mat-3-3.center)--
  (mat-4-3.center);
\draw[opacity=0.4,myarrow=0.33] 
  (mat-4-3.center) --
  (mat-4-4.center);
\draw[opacity=0.4,myarrow=0.37] 
  (mat-4-4.center) --
  (mat-2-4.center);
\draw[opacity=0.4,myarrow=0.37] 
  (mat-2-4.center) --
  (mat-2-2.center);
\draw[opacity=0.4,myarrow=0.17] 
  (mat-2-2.center) --
  (mat-5-2.center);
\draw[opacity=0.4,myarrow=0.17] 
  (mat-5-2.center) --
  (mat-5-5.center);
\draw[opacity=0.4,myarrow=0.14] 
  (mat-5-5.center) --
  (mat-1-5.center);
\draw[opacity=0.4,myarrow=0.14] 
  (mat-1-5.center) --
  (mat-1-1.center);
\draw[opacity=0.4,myarrow=0.1] 
  (mat-1-1.center) --
  ([yshift=-25pt]mat-5-1.center);
\end{pgfonlayer}
\end{scope}
\end{tikzpicture}
\]

\printbibliography

\end{document}